\documentclass{AIMS}
\usepackage[english]{babel}
\usepackage{amsmath,amsfonts,amssymb,amsopn,amscd}
\usepackage{bbm,wasysym,stmaryrd}
\usepackage[colorlinks=true]{hyperref}

\newtheorem{theorem}{Theorem}[section]

\newtheorem{lemma}[theorem]{Lemma}
\newtheorem{proposition}{Proposition}

\theoremstyle{definition}
\newtheorem{definition}[theorem]{Definition}
\newtheorem{remark}{Remark}

\newcommand{\RR}{\mathbbm{R}}
\newcommand{\R}{\mathbbm{R}}
\newcommand{\LL}{\mathbbm{L}}

\newcommand{\derpart}[2]{ \frac{\partial #1}{\partial #2} }	
\newcommand{\m}{\mathcal}
\newcommand{\vect}{\overrightarrow}
\newcommand{\interieur}{\stackrel{\mbox{\tiny o}}}
\newcommand{\inv}{\frac{1}}
\newcommand{\fom}{\m{F}(\interieur{\Omega_h})}
\newcommand{\Om}{\interieur{\Omega_h}}
\newcommand{\f}{\varphi}

\newcommand{\w}{\Omega}
\newcommand{\wt}{\widetilde{\Omega}}

\newcommand{\ut}{\tilde{u}}
\newcommand{\Lop}{\mathcal{L}}

\hypersetup{urlcolor=blue, citecolor=red}

\title[Shape-Controllability of Evolution Equations ]{Controllability of the heat and wave equations and their finite difference approximations by the shape of the domain}
\author[Jonathan Touboul]{}



 \email{jonathan.touboul@college-de-france.fr}

\begin{document}
\maketitle
\centerline{\scshape Jonathan Touboul }
\medskip
{\footnotesize
 \centerline{Mathematical Neuroscience Laboratory, CIRB-Coll\`ege de France and BANG Laboratory, INRIA Paris-Rocquencourt.}
   \centerline{11, place Marcelin Berthelot}
   \centerline{ 75005 Paris, FRANCE}
} 

\bigskip

 \centerline{(Communicated by the associate editor name)}

\begin{abstract}
	 In this article we study a controllability problem for a parabolic and a hyperbolic partial differential equations in which the control is the shape of the domain where the equation holds. The quantity to be controlled is the trace of the solution into an open subdomain and at a given time, when the right hand side source term is known. 
	The mapping that associates this trace to the shape of the domain is nonlinear. We show (i) an approximate controllability property for the linearized parabolic problem and (ii) an exact local controllability property for the linearized and the nonlinear equations in the hyperbolic case. We then address the same questions in the context of a finite difference spatial semi-discretization in both the parabolic and hyperbolic problems. In this discretized case again we prove a local controllability result for the parabolic problem, and an exact controllability for the hyperbolic case, applying a local surjectivity theorem together with a unique continuation property of the underlying adjoint discrete system.
\end{abstract}

\section*{Introduction}
The problem of characterizing the shape of a domain where a certain dynamical phenomenon is partially observable is a model for a wide class of applications. A typical example is given by the identification of the shape of an hydrocarbon or water reservoir. Techniques used by geologists to tackle this problem include sending shock waves into the ground or into a drill hole and measure the reflected waves. The question we may address here is whether one can infer the shape of the full reservoir by these partial and local informations. Chenais and Zuazua \cite{CZ} addressed this problem for the elliptic equation case dealing with the Laplace equation with Dirichlet boundary conditions. They showed that the linearized problem admits an approximate controllability property, and the finite-differences discretization presents a local controllability property. 

In the present manuscript, we extend their results by addressing the same questions in a dynamical setting. In details, we shall consider in this manuscript evolution partial differential equations of parabolic and hyperbolic type equation on an open set $\Omega \subset \RR^n$, given an external source term. The domain $\Omega$ is assumed to be only partially known. It is potentially allowed to evolve with time, in which case it will be denoted $\Omega(t)$, and is assumed to contain for all times a fixed simply connected open subset $\omega \subset \RR^n$. More precisely, we assume that the closure of $\omega$ is a subset of $\Omega(t)$ for all times, a property we shall denote $\omega \Subset \Omega$. We assume that the solution of the partial differential equation restricted to $\omega$ at time $T$ is known, either because it is accessible or observable by means of measurement. The problem we address is, loosely speaking, to recover the shape $\Omega(t)$ from the knowledge of the external forcing term and the restriction of the solution to the subdomain $\omega$ at time $T$. This is the question we refer to as the 
\emph{controllability problem}. It differs from \cite{CZ} in that we consider evolution equations and allow the domain to dynamically evolve in time. These distinctions yield increased complexity of the functional setting that necessitate to partially modify the result obtained by Chenais and Zuazua.   

In this article, we only address the case of the heat and wave equations as paradigmatic examples of, respectively, parabolic and hyperbolic equations and will restrict the study to the case of Dirichlet boundary conditions, although the problem might arise for other differential operators and different boundary conditions and the methods used here may apply in most regular cases. 

Returning to the mathematical formulation of the problem, we are given:
\begin{itemize}
	\item $\Lop$ a given parabolic or hyperbolic differential operator, i.e. either the heat or the wave operator,
	\item a source-term $f\in \LL^2(\R^+\times \R^n)$,
	\item $\omega$ an open bounded subset of $\R^n$ assumed regular,
	\item $y_d \in H^1(\omega)$ the observed solution of the PDE $\Lop y =f$ at time $T$ restricted to $\omega$ .
\end{itemize}
and aim at proving existence and uniqueness of a bounded, possibly time-varying, open set $\Omega(t)\subset \R^n$ for $t\in [0,T]$ with $\omega \Subset \Omega(t)$ for all $t$ and such that the solution $y_{\Omega(\cdot)}$ of the equation:
\[\begin{cases}
	\Lop y_{\Omega} = f & \text{on } \;\Omega(t)\\
	y_{\Omega}(t,x) = 0 & x\in \partial \Omega(t)
\end{cases}\]
satisfies $y_{\Omega}(T,\cdot)\vert_{\omega}=y_d$. 

It is therefore a shape identification problem, which can be seen as a controllability problem in the sense that the domain $\Omega$ has to be determined so that $y_{\Omega}(T,\cdot)\vert_{\omega}=y_d$ holds. Problems of this type, in different settings, have been addressed using a variety of methods. Some authors have used optimization methods (e.g. \cite{C,chenais:75, cea:81, MS, rousselet:81, simon:80,zolesio:81}), writing the problem in the form:
\[\inf\limits_{\Omega \in \m{U}} \|y_{\Omega}-y_d\| \] 
where $\m{U}$ is the set of domains we take into account. Under suitable conditions on the set of admissible domains $\m U$, they show that the existence of a minimizer can be guaranteed. However, this existence result does not ensures that $y_{\Omega}\vert_{\omega}= y_d$ does actually hold for a certain choice of the domain $\Omega$, nor does it evaluate the minimal distance between $y_{\Omega}\vert_{\omega}$ and $y_d$. Therefore, optimization techniques will not solve the controllability problem under consideration.

It is important to note the fact that solutions to the controllability problem do not necessarily exist, and if they exist, these are not necessarily unique. In order for our controllability problem to be well posed, we will be interested in a local controllability property. More precisely, we consider a reference domain $\Omega_0$ for which the partial differential equation has a solution $y_0$. The question we ask is the following: given $y_d$ a function close from $y_0(T)\vert_{\omega}$ in a suitable topology, is there an open set $\Omega^*(t)$ close from $\Omega$ (in a sense that will be defined in the sequel) such that $y_{\Omega^*}(T)\vert_{\omega} = y_d$. This problem is referred to as the \emph{local controllability problem}. It may also happen that we can approach $y_d$ arbitrarily close but cannot find an optimal open set $\Omega^*(t)$ such that $y_{\Omega^*}(T)\vert_{\omega} = y_d$. We will hence consider:
\begin{itemize}
	\item the \emph{exact} controllability properties, where there exists a domain $\Omega^*(t)$ such that $y_{\Omega^*}(T,\cdot)\vert_{\omega} = y_d$
	\item and the \emph{approximate} controllability property, where for any $\varepsilon>0$ there exists a domain $\Omega_{\varepsilon}(t)$ such that $\big \Vert y_{\Omega_{\varepsilon}}(T,\cdot)\vert_{\omega} - y_d\vert_{\omega}\Vert \leq \varepsilon$ in a suitable space. 
\end{itemize}
There is an extensive literature on exact and approximate controllability problems for partial differential equations (see e.g. \cite{L,Z1, Z2}), but very little has been done for controls with respect to the shape of the domain.

As discussed in detail in~\cite{CZ}, the problem we consider presents various technical difficulties preventing the use of existing methods. The main problem is the nonlinear nature of the mapping $\Omega\mapsto y_{\Omega}$, differing from existing results in nonlinear controllability problems often valid only for ``mild'' nonlinear perturbations, and which is not the case here since the trace of the solution depends on the shape of the domain in a genuinely nonlinear way. 

In order to achieve our program of controlling such equations, it is natural to start by linearizing the problem, with the aim of applying the Inverse Function Theorem (IFT). This approach, developed in section~\ref{sect:continuous}, allows proving that the linearized control problem presents an \emph{approximate local} controllability property in the parabolic case and a \emph{exact local } controllability property in the hyperbolic case. The proof of this property consists in applying a duality argument together with a unique continuation result for the solutions of the adjoint system, which is obtained as a consequence of Holmgren's uniqueness theorem.

However, even if we are able to prove the linearized problem is approximately or exactly controllable, this does not allow to conclude a controllability result about the initial nonlinear system. The main limitations are related to the inherent complexity of the spaces in which the nonlinear problem holds. These limitations are released when considering spatially discretized versions of the problem, holding in simpler finite-dimensional spaces. This is what we show in section \ref{sect:semi-discrete} in the context of semi-discrete finite-difference approximations of the initial problem. By reducing the problem into an ordinary differential equation in a finite-dimensional space, we will prove local controllability property for the semi-discrete parabolic problem and exact controllability result in the hyperbolic case. In the context of finite-dimensional systems, approximate and exact controllability are equivalent notions, and the problem may be reduced to an unique continuation issue for the adjoint system for the linearized system. The price to pay is that classical tools used for unique continuation used in the domain of PDEs (like Holmgren Uniqueness Theorem or Carleman inequalities) do not seem to apply in spatially discrete systems. Thus, the first task that we undertake is to prove a new unique continuation result holding in our setting based on propagation properties the discrete scheme naturally induces. This allows to demonstrate that the linearized model is exactly controllable, which implies a local controllability result in virtue of the IFT. Note that these results are only concerned with local controllability results for fixed discretization meshes. A very interesting problem would be to address the convergence of the control shapes as the mesh-size $h$ tends to zero. This is a problem of primary importance but its analysis is out of reach for the techniques developed here and not in the scope of this paper. At this respect it is important to note that very little is known about the convergence of controls in the context of the controllability of numerical approximations of PDEs too. For instance, in the case of controllability of the wave equation where the control is on the boundary condition, it is known that controls do not necessarily converge as the mesh size tends to zero because of high frequency spurious oscillations (see~\cite{Z4}). However, in the context of the heat equation, at least in one dimension, the controls driving solutions to rest do converge as the mesh size tends to zero~\cite{Z6}. Similar results hold in the context of homogenization for wave and heat equations with rapidly oscillating coefficients~\cite{Z3}. It is highly nontrivial to extend these techniques when the control is the shape of the domain. 

%

\section{Controllability of the heat and wave equations}
\label{sect:continuous}
We consider here the problem of the controllability of the heat and wave PDEs as a function of the shape of the domain. 

\subsection{Controllability of the heat equation}


Let us start by introducing the mathematical framework used throughout this section. We classically denote by $W^{k,\infty}(\mathbb{R}^n,\mathbb{R}^m)$ the set of functions $k$ times differentiable of $\mathbb{R}^n$ (in the sense of distributions), taking values in $\mathbb{R}^m$, with all differentials in $L^{\infty}(\RR^n,\RR^m)$. We consider the classical heat equation with a source term $f \in  L^2(\mathbb{R}^+, \mathbb{R}^n)$ holding on an open set $\Omega_0$ with smooth boundary $\partial\Omega_0 \in W^{2,\infty}$: 

\begin{equation}\label{eq:BaseHeat}
 \begin{cases}
        \partial_t y_0(t,x) - \Delta y_0(t,x)  = f(t,x)  & t>0 \textrm{ and } x\in \Omega_0 \\
        y_0(0,x)= 0 & x \in \Omega_0 \\
        y_0(t,x) = 0 & t>0 \textrm{ and } x \in \partial\Omega_0 .
 \end{cases}
\end{equation}
This equation defines a unique solution $y_0$. We are interested in possible values of $y(T)\vert_{\omega}$ solutions of the heat equations holding in $\Omega(t)$ a ``perturbation'' of $\Omega_0$. We work in the standard setting for differentiation with respect to the domain (see for instance \cite{MS, rousselet:81, S,simon:80, zolesio:81}). The admissible perturbed open sets are ``small dynamical perturbations'' of $\Omega_0$. In details, we consider deformation functions $\varphi:\R^+\times \R^n\mapsto \R^n$ that are $\LL^2([0,T], W^{2,\infty}(\mathbb{R}^n,\mathbb{R}^n))$. We moreover consider transformations such that for all $t \in [0,T]$, the maps $(id + \varphi(t))$ and $(id +\varphi(t))^{-1}$ are homemorphisms of $W^{2,\infty}(\mathbb{R}^n,\mathbb{R}^n)$\footnote{This is always the case if for all $t\in [0,T]$, $\varphi(t)$ is close from $0$ in $W^{2,\infty}(\mathbb{R}^n,\mathbb{R}^n)$.}. Our choice will be, within this set, to consider $\varphi$ such that $\Vert \varphi\Vert \leq \varepsilon$ for some $\varepsilon>0$. This set is noted $\m W$.

For $\varphi \in \m{W}$, we define
\[\Omega_{\varphi}(t) = (id + \varphi(t)) \, (\Omega_0) = \{ x + \varphi(t)(x) ; x \in\Omega_0\}.\] 
The perturbed equation is:
\begin{equation}\label{eq:PerturbedHeat}
 \begin{cases}
        \partial_t y_{\varphi}(t,x) - \Delta y_{\varphi}(t,x)  = f(t,x)  & t>0 \textrm{ and } x \in \Omega_{\varphi}(t)\\
        y_{\varphi}(0,x)= 0 & x \in \Omega_{\varphi}(0) \\
        y_{\varphi}(t,x) = 0 & t>0 \textrm{ and } x\in\partial\Omega_{\varphi}(t)
 \end{cases}
\end{equation}
and function associating the trace of the solution on $\omega$ is noted $\Lambda$:
\begin{equation}
\Lambda:= 
\left \{
\begin{array}{ll}
    \m{W} &\mapsto  H^1(\omega)\\
    \f    &\mapsto   y_{\f}(T,\cdot)\vert_{\omega}.
\end{array}
\right .
\end{equation}
The range of $\Lambda$: $\m{R}(\Lambda) = \Lambda(\m{W})=\{ y_{\f}(T,\cdot)\vert_{\omega} \, ; \f \,\in \m{W}\}$ constitute the set of accessible states at time $T$.


\subsubsection{Existence and uniqueness of solutions for the perturbed problem}
In this section we use the variational formulation and show weak existence and uniqueness of the solution of the perturbed system. 

\begin{lemma}\label{theo:HeatEquivalentProblem}
 The perturbed problem \eqref{eq:PerturbedHeat} is equivalent to the following partial differential equation on $\overline{y_{\f}}(t,x) = y_{\varphi}\Big(t,x+\varphi(t)(x)\Big)$ on $\R^+\times \Omega_0$:
\begin{equation}\label{eq:equivalentProblem}
  \begin{cases}
     \displaystyle{\partial_t \overline{y_{\f}} - \frac{1}{\vert \textrm{det}(id + \nabla \varphi (t))\vert} \textrm{div}(B(\varphi)\nabla \overline{y_{\f}}) = f \circ (id + \varphi)} & t>0 \textrm{ and } x \in \Omega_0\\
     \overline{y_{\f}}(0,x) = 0 & x\in\Omega_0\\
     \overline{y_{\f}}(t,x) = 0 & t>0 \text{ and } x\in{\partial\Omega_0}
  \end{cases}
\end{equation}
with $f\circ (id+\varphi) (t,x)=f(t,x+\varphi(t)(x))$ and 
\begin{equation}\label{eq:BFunction}
 B(\varphi) = \Big \vert\textrm{det} \big (id + \nabla \varphi (t)\big ) \Big\vert\; \left ( \Big [ \big ( \nabla(id + \varphi(t) ) \big)^*\Big]^{-1}\right )^* \; \Big[\big (\nabla(id + \varphi(t))\big)^*\Big]^{-1}.
\end{equation}
\end{lemma}

\begin{proof}
We use the variational formulation of equations \eqref{eq:PerturbedHeat} and show that it is equivalent to the variational formulation of the problem \eqref{eq:equivalentProblem}. 
Let $\phi \in H^1 \Big ( \mathbb{R}^+, H^1_0(\Omega_{\varphi}(t)) \Big)$ a test function. We have:
\begin{align*}
    \int_{\mathbb{R}^+ \times
    \Omega_{\varphi}(t)}f \phi dt\, dx &= \int_{\mathbb{R}^+ \times \Omega_{\varphi}(t)} (\partial_t
    y_{\varphi}\phi - \Delta y_{\varphi} \phi ) \,dt\, dx \\
    &= \int_{\mathbb{R}^+ \times \Omega_{\varphi}(t)} (\partial_t
    y_{\varphi} \phi + \nabla y_{\varphi} \nabla \phi)\, dt\, dx.
\end{align*}
We introduce the function $\Psi$ defined on $\mathbb{R}^+ \times \Omega_0$ by $\Psi(t,x)= \phi(t, x+\varphi(t)(x))$, and change variables in the above weak formulation defining $r\in\Omega_0$ such that $x=r+\varphi(r)$. The determinant of the Jacobian matrix of the change of variable is given by $\textrm{det}(id + \nabla \varphi)$, yielding:
\begin{multline}
    \int_{\mathbb{R}^+ \times \Omega_0} \left \{ \partial_t
    \overline{y_{\f}} \Psi - (\nabla y_{\varphi}) \circ (id + \varphi) \cdot
    (\nabla \phi)\circ (id + \varphi) \right \} \left \vert \textrm{det}(id + \nabla \varphi (t) ) \right \vert dt\, dx \\
    = \int_{\mathbb{R}^+ \times \Omega_0} f\circ (id + \varphi) \Psi \left \vert\textrm{det}(id + \nabla \varphi (t)) \right \vert \, dt\, dx
\end{multline}
We therefore need to express $(\nabla \Phi) \circ (id + \varphi)$ as a function of $\nabla \Big (\Phi \circ (id + \varphi) \Big)$. Let $H=\Phi \circ (id + \varphi)$ and $T=(id + \varphi)$. We clearly have $\nabla H = (\nabla T)^*\,\nabla \Phi$
%
where the star denotes the adjoint operator, yielding:
%
%
%
\begin{multline}
    \int_{\mathbb{R}^+ \times \Omega_0} \partial_t \overline{y_{\f}} \Psi \,-\, \Big[ \big (\nabla(id + \varphi(t)) \big )^*\Big ]^{-1} \nabla \overline{y_{\f}} \Big [ \big (\nabla(id + \varphi(t))\big)^*\Big]^{-1} \nabla \Psi \left \vert\textrm{det} (id + \nabla \varphi (t)) \right \vert \, dt\, dx \\ 
    = \int_{\R^+ \times \Omega_0} f\circ (id + \varphi) \Psi \left \vert\textrm{det} (id + \nabla \varphi (t)) \right \vert \, dt\, dx.
\end{multline}
Using the operator $B$ introduced in equation \eqref{eq:BFunction}, it is easy to see that the initial variational problem is equivalent to the variational problem:
\begin{multline}
\int_{\R^+ \times \Omega_0} \left \{ \partial_t \overline{y_{\f}} \, \left |\textrm{det} (id + \nabla \varphi (t)) \right| - \textrm{div}\big(B(\varphi)\big)\nabla \overline{y_{\f}}\right\} \Psi dt\, dx\\
     = \int_{\R^+ \times \Omega_0} f \circ (id + \varphi) \Psi \left | \textrm{det} \big(id + \nabla \varphi (t)\big )\right |dt\, dx
\end{multline}
whose solutions define the weak solutions of the partial differential equation \eqref{eq:equivalentProblem}. 
\end{proof}

This formulation has the interest to hold on a fixed domain. Based on this formulation, we show the well-posedness of the perturbed problem.
\begin{proposition}\label{theo:existenceUniquenessHeat}
 The PDE \eqref{eq:PerturbedHeat} has a unique solution in $L^2\Big(0,T; H^1_0(\Omega_0)\Big) \cap C\Big([0,T],L^2(\Omega_0)\Big)$.
\end{proposition}
%
%
%
%
%
%
\begin{proof}
Let $H = L^2(\Omega_0)$ and $V = H^1_0(\Omega_0)$, and recall that the dual space of $V$ is, by definition, $V'= H^{-1}(\Omega_0)$. Let
\[\begin{cases}
	f & \in L^2(0,T, V')\\ 
	u_0 & \in L^2(\Omega_0) = H
\end{cases}
\]
A classical result due to Lions (see e.g.~\cite[theorem X.9]{B}) ensures that, provided that there exists a function  $a(t;u,v)$ measurable as a function of time, bilinear in $(u,v)$ for all $t\in [0,T]$ and such that there exist constants $M, \alpha$ and $C$ such that:
\[\begin{cases}
	    |a(t;u,v)| &\leq M \|u\| \|v\| \\
	    a(t;u,v) &\geq \alpha \|v\|_V^2 - C |v|_H^2
\end{cases}
	\]
then there exists a unique function $u$ such that: 
\begin{itemize}
    \item $u \in L^2(0,T; H^1_0(\Omega_0)) \cap C([0,T],L^2(\Omega_0))$
    \item $\partial_t u \in L^2(0,T;V')$
    \item $
        \left \{
            \begin{array}{cccl}
                <\partial_t u, v> + a(t;u(t), v) &=& <f(t),v>
                & \textrm{for a.e. } t\in (0,T), \forall v \in H^1_0(\Omega_0)\\
                u(0) &=& u_0 & \textrm{ on   }\Omega_0.
            \end{array}
        \right . $
\end{itemize}

This theorem readily applies to our case. Indeed, the bilinear form on $L^2(\Omega_0)$ given by: 
\[\langle f,g\rangle = \int_{\Omega_0} f(x) g(x) |\hbox{det}(id + \nabla \varphi)|  dx .\]
is positive definite because of the invertibility of $id+\varphi(t)$ for $\f\in\m{W}$, hence defines a suitable scalar product on $L^2(\Omega_0)$. Let us introduce the function: 
\[ a(t;u,v) = \int_{\Omega_0} B(\varphi) \nabla v \nabla \left ( \frac{u}{|\hbox{det}(id + \nabla \varphi)|} \right ) dx,\] 
which is well defined since we assumed $\f \in W^{2,\infty}(\R^n,\R^n)$. The function $a$ is clearly measurable in $t$. Moreover, for $\varphi$ close enough from $0$ in $\m W$, we have $B(\varphi)$ close from identity, and hence $\left | \; \|B(\varphi)\|-1 \; \right|\leq 1/2$ for $\varphi$ small enough. We also have $|\textrm{det}(id + \nabla \varphi)| \to 1$ when $\varphi \to 0$ in $\m{W}$. 

Therefore, by sufficiently restricting the space $\m W$, the following two inequalities simultaneously hold: 
\begin{itemize}
	\item $\vert \;\vert \hbox{det}(id + \nabla\varphi)\vert-1\vert \leq 1/2$, and hence $|\hbox{det}(id + \nabla \varphi)| \geq 1/2$
	\item $\|B(\varphi)\| \leq 3/2$.
\end{itemize} 
Moreover the function $\frac{u}{|\hbox{det}(id + \nabla \varphi)|}$ belongs to $H^1_0(\Omega_0)$, implying:
\begin{eqnarray*}
    |a(t;u,v)| &\leq& \int_{\Omega_0} \left | B(\varphi) \nabla
    v \nabla \left ( \frac{u}{|\hbox{det}(id + \nabla \varphi)|} \right ) \right | dx\\
    {} & \leq& 3/2 \|\nabla v\| \; \left \|\nabla \left (\frac{u}{|\hbox{det}(id + \nabla
    \varphi ) | } \right) \right \|\\
    {} & \leq&  3\| v\|_{H^1_0} \; \|u\|_{H^1_0}
\end{eqnarray*}
and eventually, 
\begin{eqnarray}
    \nonumber a(t;v,v) &=& \int_{\Omega_0} B(\varphi) \nabla v \nabla \left (\frac{v}{|\hbox{det}(id + \nabla \varphi)|}\right )\\
    \nonumber {}&\geq & \frac 1 2 \int_{\Omega_0} \nabla v \nabla \left (\frac{v}{|\hbox{det}(id + \nabla \varphi)|}\right )\\
    \nonumber {} &\geq & \frac 1 3 \|v\|^2_{H^1_0(\Omega_0)}.
\end{eqnarray}
Lions theorem hence applies, ensuring the existence and uniqueness of a solution to the equations \eqref{eq:equivalentProblem} and this solution belongs to the space $L^2(0,T; H^1_0(\Omega_0)) \cap C([0,T],L^2(\Omega_0))$. This directly implies existence and uniqueness of a solution in $L^2(0,T; H^1_0(\Omega_{\varphi})) \cap C([0,T],L^2(\Omega_{\varphi}))$ for equations~\eqref{eq:PerturbedHeat}. 
\end{proof}



\subsubsection{Linearized problem}
Now that we have shown that the nonlinear problem was well-posed, we introduce and analyze the linearized problem. The linearization is performed with respect to the parameter $\varphi$. There are different notions of differentials with respect to the shape of the domain. In this article we are interested in the Lagrangian differential, denoted $y_{\varphi}'$, defined as the unique function (when it exists) satisfying the property:
\[\forall \widetilde{\Omega} \Subset \Omega_0 \;\;
y_{\varphi}|_{\widetilde{\Omega}}=y_{0}|_{\widetilde{\Omega}}+y_{\varphi}'|_{\widetilde{\Omega}}
+o(\varphi).\]

\begin{proposition}
The Lagrangian shape differential $y_{\varphi}'$ satisfies the following equation:
\begin{equation}\label{eq:HeatDerLag}
        \begin{cases}
            \partial_t y_{\varphi}'(t,x) - \Delta y_{\varphi}'(t,x) = 0 &  t>0 \text{ and } x \in \Omega_0\\
            y_{\varphi}'(0,x)=0 & x \in \Omega_0\\
            y_{\varphi}'(t,x) = -\varphi(t).n \; \derpart{y_0}{n}& t>0 \text{ and } x \in \partial \Omega_0.
        \end{cases}
\end{equation}
where the dot denotes the scalar product in $\R^n$.
\end{proposition}

\begin{proof}
Let $\widetilde{\Omega} \Subset \Omega_0$. We have in this open set the evolution equation
\[ \partial_t y_{\varphi} - \Delta y_{\varphi}  = f  \;\; \forall (t,x) \in \RR^+ \times \widetilde{\Omega}.\] 
Since $y_0$ satisfies the equation:
\begin{eqnarray*}
   \partial_t y_{0} - \Delta y_{0} & = f & \textrm{ on } \; \RR^+ \times \widetilde{\Omega},
\end{eqnarray*}
we necessarily have the following evolution equation for $y_{\varphi}'$ on $\RR^+ \times \widetilde{\Omega}$:
\begin{equation}\label{eq:HeatDerLagrange}
  \partial_t y_{\varphi}' - \Delta y_{\varphi}'  = 0  \;\;\; \textrm{on }\; \RR^+ \times \widetilde{\Omega}.
\end{equation}
The evolution equation of the Lagrangian shape differential being derived, we now identify the related boundary conditions. First of all, it is clear that 
\[y_{\varphi}'(t=0,x)=0 \qquad \textrm{ on } \widetilde{\Omega}.\] 
Let us now compute $y_{\varphi}|_{\partial\Omega_0}$. Considering, for some $\psi \in C^{\infty}(\mathbb{R}^n)$ fixed, the integral $\int_{\partial\Omega}y_{\Omega}(t) \psi ds$ and formally differentiating with respect to $\Omega$, we obtain the following formula (this calculation is rather classical, see for instance \cite[6.28]{A}):
\begin{eqnarray*}
    0 &=& \left . \frac{\partial}{\partial \Omega}\right
    |_{\Omega =
    \Omega_0} \int_{\partial \Omega} y_{\Omega}(t)\psi ds\\
    &=& \int_{\partial  \Omega_0} y_{\varphi}'(t) \psi
    ds + \int_{\partial  \Omega_0} \varphi(t).n \frac{\partial (y_0 \psi)}{\partial n}
    ds\\
    &=& \int_{\partial  \Omega_0} y_{\varphi}'(t) \psi
    ds + \int_{\partial  \Omega_0} \varphi(t).n (\frac{\partial y_0 }{\partial
    n}\psi +y_0 \frac{\partial \psi}{\partial n})
    ds\\
    &=& \int_{\partial  \Omega_0} y_{\varphi}'(t) \psi
    ds + \int_{\partial  \Omega_0} \varphi(t).n \frac{\partial y_0 }{\partial n}\psi
    ds
\end{eqnarray*}
\noindent because $y_0|_{\partial  \Omega_0}\equiv 0$. Hence we have:

\[ \int_{\partial \Omega_0} (y_{\varphi}'(t) + \varphi(t).n \frac{\partial y_0}{\partial n})\psi ds = 0 \;\;\; \forall \psi \in C^{\infty}(\RR ^2).\]
The Lagrangian shape differential therefore satisfies the non homogeneous Dirichlet boundary condition:
\[ y_{\varphi}'(t,x) = -\varphi(t)\cdot n \frac{\partial y_0}{\partial n} (t,x)\qquad \textrm{ for } (t,x)\in{\R^+\times\partial\Omega_0}\]
yielding our equation \eqref{eq:HeatDerLag}.
\end{proof}


\subsubsection {Approximate controllability of the linearized problem}

If $d\Lambda(0)$ was an isomorphism of $W^{2,\infty}(\RR^n, \RR^n)$ on $H^1(\RR^n)$, then the local inversion theorem would readily imply exact local controllability. However we will see that this property does not hold true, and a weaker property will be proved: we show that the linearized problem is approximately controllable. This property is demonstrated using Holmgren's theorem~\cite{H,Holmgren:1901} giving a uniqueness result for PDEs with real analytic coefficients, and will be used here to show a propagation of zeros property from a non-characteristic surface. 
%
%

\begin{lemma}\label{prop:HolmgrenHeat}
Let $\Omega$ an open subset of $\RR^n$, with regular boundary $\partial\Omega$ (e.g. $W^{2,\infty}$ as assumed for $\Omega_0$). Let $\gamma$ a non-empty open subset of $\partial\Omega$. Any solution of the equations:
\begin{equation}\label{eq:EDPHeatu}
  \begin{cases}
    \partial_t u(t,x) - \Delta u(t,x) = 0 & t>0, x\in \Omega\\
    u(0,x)=0  & x \in \Omega\\
    u (t,x) = 0 & t>0, x\in \partial\Omega_0\\
    \frac{\partial u}{\partial n} =0 & \forall t>0, x\in \gamma
  \end{cases}
\end{equation}
is null on $\R^+\times \Omega$.
\end{lemma}

\begin{proof}
Let $D=(\partial_t,\partial_{x_1} ,\ldots, \partial_{x_n})$ and $P(D)= -\partial_t - \Delta$ the heat differential operator. The associated characteristic polynomial is simply $P(T,X) = -T - |X|^2,$ and its principal part is $P_2(T,X)=-|X|^2=0$. The solution of this later equation is $X=0$, and its direction $(1,0,\ldots,0)$, hence the characteristic surfaces are the hyperplanes $T=\textrm{ constant }$. 
The variational formulation satisfied by the function $u \in C(0,T;H^1_0(\Omega))$ solution of \eqref{eq:EDPHeatu} is given by:
\begin{equation}\label{eq:HeatuVariationnal}
    \forall v \in H^1(\Omega),
    \int_{\Omega} (\partial_t u \, v + \nabla u \, \nabla v) \, dx -
    \int_{\partial\Omega} \frac{\partial u}{\partial n} v \, ds = 0.
    \end{equation}
Let now  $\tilde{u}$ be a continuation of $u$ on an open subset $\widetilde{\Omega}=\w \cup \Gamma$ where  $\Gamma$ is an open set whose intersection with $\Omega$ is equal to $\gamma$.

Let $\tilde{u} \in H^1_0(\wt)$ be the function defined by $\ut(t,x) = u(t,x) \mathbbm{1}_{x\in\w}$. We have the following relations:
    \begin{equation*}
     \begin{cases}
            \partial_t \ut(t,x) = 0 & \textrm{ on }\R^+\times {\Gamma}\\
            \nabla \ut = 0 & \textrm{ on }\R^+\times {\Gamma}\\
            \frac{\partial \ut}{\partial n} (t,x)= 0 &\textrm{ on } \gamma
     \end{cases}
    \end{equation*}

    \noindent and hence we have  :
    \begin{equation}\label{eq:utVariational}
        \forall v \in H^1(\Gamma),
        \int_{\Gamma} \partial_t \ut \, v + \nabla \ut \, \nabla v -
        \int_{\partial \Gamma} \frac{\partial \ut}{\partial n} v ds = 0.
    \end{equation}

    \noindent Now, using \eqref{eq:HeatuVariationnal} and \eqref{eq:utVariational}, the condition $\frac{\partial u}{\partial n}\vert_{\gamma}\equiv0$ and the fact that $\Gamma \cup \partial \Omega = \gamma$, we obtain that  $\ut$ satisfies the following variational problem:
    \begin{equation}\label{fvut}
        \forall v \in H^1(\wt),
        \int_{\wt} \partial_t \ut \, v + \nabla \ut \, \nabla v -
        \int_{\partial \wt} \frac{\partial \ut}{\partial n} v ds = 0.
    \end{equation}
    
    If both $\wt$ and $\Gamma$ were convex sets, Holmgren's theorem readily implies that $\ut\equiv 0$. Indeed, $\R^+\times\Gamma\subset \R^+\times\wt$, $P(D)$ is a differential operator with constant coefficients and every plane which is characteristic with respect to $P$ intersecting $\R^+\times\wt$ intersects $\R^+\times\Gamma$ as well (since these are the planes of constant $t$). Under these conditions, Holgren's theorem~\cite[Theorem 5.3.3]{H} implies that any solution of  $P(D)u=0$ on $\wt$ vanishing on $\R^+\times\Gamma$ also vanishes on $\R^+\times\wt$, and in particular $u\equiv 0$.

In our framework, it can occur that the sets $\Gamma$ and $\wt$ are not convex. However, it is always possible to describe $\wt$ as the union of open balls (because of the regularity and connectedness of $\w$). These discs are convex, and the intersection between two discs is also convex, and the above argument applied on each element of this decomposition yield the desired result. 
    
\end{proof}

Thanks to these results, we can prove the density result previously announced for the continuous problem:

\begin{theorem}\label{theo:densite_continu}
    Assume that there exist a non-empty subset $\gamma$ of $\partial\omega$ on which 
    \[\forall (t,x) \in [0,T]\times \gamma ,\;\;\; \frac{\partial y_0}{\partial n} (t,x)\neq 0\]
    \noindent Then $\mathcal{R}=\{y_{\varphi}'(T)\vert_{\omega}; \varphi  \in \m{W}\}$ is dense in $L^2(\omega)$.
\end{theorem}

\begin{proof}
We prove that $\mathcal{R}^{\bot}=\{0\}$. Indeed, any $g \in L^2(\Omega_0)$ belonging to $\mathcal{R}^{\bot}$ is such that:
$$\int_{\omega} g \; h =0 \textrm{, } \forall h \in \mathcal{R}.$$ 
By definition, $\mathcal{R} = \{y_{\varphi}'(T)\vert_{\omega}; \varphi \in \m{W}\}$ so the later condition is equivalent to :
\begin{equation}\label{eq:surjective}
\forall \varphi \in \m{W} \int_{\Omega_0} g y_{\varphi}'(T) = 0.
\end{equation}
Let us define the adjoint state $\phi \in H^1(\Omega)$ associated to $g$:
\begin{equation*}
 \begin{cases}
    -\partial_t \phi(t,x) - \Delta \phi(t,x) &= g \otimes \delta_{t=T} \qquad t\geq 0, x\in \omega\\
    \phi(t=0,x)&=0 \qquad x\in \omega\\
    \phi(t,x) &= 0 \qquad t\geq 0, x\in \partial\omega
 \end{cases}
\end{equation*}
where we denote:
\[\forall v \in C(0,T;H^1(\omega)) \qquad \langle g \otimes \delta_{t=T}, v\rangle = \int_{\omega} g v(T) dx.\] 
We now show that $\phi = 0$, that is $g=0$. Equation \eqref{eq:surjective} implies that for any $\varphi \in \m{W}$, 
$$ \langle -\partial_t \phi - \Delta \phi, y_{\varphi}'\rangle =
\int_{\omega} g y_{\varphi}'(T) =0.$$ 
Integrating by parts yields: 
\begin{align*}
    0=\int_{\mathbb{R}^+ \times \omega} (-\partial_t \phi - \Delta \phi) y_{\varphi}' &= \int_{\mathbb{R}^+ \times    \omega} \phi \partial_t y_{\varphi}' - \int_{\omega} \phi(0) y_{\varphi}'(0) dx \\
    & + \int_{\mathbb{R}^+ \times \omega}\nabla \phi \nabla y_{\varphi}'  - \int_{\mathbb{R}^+ \times \partial\omega} y_{\varphi}' \frac{\partial \phi}{\partial n}\\
    {}&= \int_{\mathbb{R}^+ \times \omega} (\partial_t y_{\varphi}' - \Delta y_{\varphi}')\phi + \int_{\mathbb{R}^+ \times \partial\omega} -y_{\varphi}'\frac{\partial \phi}{\partial n}+\phi \frac{\partial y_{\varphi}'}{\partial n}\\
    {} &= \int_{\mathbb{R}^+ \times \partial\omega} \varphi.n \frac{\partial y_0}{\partial n} \frac{\partial \phi}{\partial
    n} dt ds.
\end{align*}
We deduce that for all $t$ in $[0,T]$, we have $\frac{\partial y_0}{\partial n} \frac{\partial \phi}{\partial n}=0 $ on $\partial\omega$.\\

Since on the non-negligible subset $\gamma \in \partial\omega$ we have $\frac{\partial y_0}{\partial n} \neq 0$ $\forall t \in [0,T]$, we deduce that necessarily $\derpart{\phi}{n}(t,x)=0$ on $\R^+\times \gamma$. Lemma~\ref{prop:HolmgrenHeat} implies that $\phi(t,x)\equiv 0$ for all $(t,x) \in (0,T)\times \omega$ or $(T,\infty)\times \omega$. The equation $$-\partial_t \phi - \Delta \phi = g \otimes
    \delta_{t=T}$$ imposes a jump condition on $\phi$ at $T$: 
    $$[| \phi|] = \phi(T^+)-\phi(T^-) = g.$$
But since here, $\phi(T^-,x)=\phi(T^+,x)=0$ for every $x\in\omega$, we necessarily have $g=0$. 
We therefore conclude that the range $\mathcal{R}$ of the operator $d\Lambda(0)$ is dense in $\LL^2(\omega)$.
\end{proof}

We emphasize the fact that the control can be chosen constant in time, i.e. that the heat equation is approximately controllable by a rigid deformation of the open set $\Omega_0$. In other words, we showed the existence of approximate controls through time invariant sets $\Omega$. 

Note also that the condition of theorem \ref{theo:densite_continu}, namely the fact that the normal differential of $y_0$ is not vanishing on a non-empty subset of the boundary of $\omega$ for all times, can appear relatively strong. This condition is however necessary in order for the observation to be performed at a given time $T$. If we are interested in the trace of the solutions on $\omega$ depending on time, $y_\f(t,x)\vert_{x\in\omega}$ and $\f\in W^{2,\infty}(\R^n,\R^n)$ (i.e. $\f$ does not depend on time), then the controllability property can be proved, using the same techniques, under the weaker assumption that $\exists \delta>0$ such that $\forall t \in [0,\delta]$  we have the non-degeneracy condition $\frac{\partial y_0}{\partial n} \neq 0$. 


We hence proved that the linearized heat equation is approximately controllable, i.e. that it has a dense range on $H^1_0(\Omega_0)$, and no exact controllability property was proved. And for good reason: the regularity introduced by the parabolic form of the equation prevents from such an exact controllability property to hold. This weaker form of controllability prevents us from using this result to address the controllability of the nonlinear problem.

%
The hyperbolic operators do not enjoy the same regularization properties as parabolic operators. We now treat the problem in the hyperbolic setting and show that, in contrast, we obtain an exact controllability property.


\subsection{Controllability of the wave equation}

We now address the same problem in the case of the wave equation. The setting and intermediate results are similar as those of the parabolic case, and will be presented in less detail. It is however important to note that the main argument ensuring controllability is completely distinct: Holmgren's uniqueness theorem was used in the parabolic case to show approximate controllability, and here we will use results on the controllability of the wave equations with respect to boundary conditions. 

We consider admissible domains $\mathcal{U}$ as dynamical perturbations of the original domain $\Omega_{\f}(t)=(id+\f(t))(\Omega_0)$ where $\f(t)$ are $\LL^2((0,T),W^{2,\infty}(\R^n,\R^n))$ such that for all $t\in[0,T]$, $\f(t)$ belongs to be the ball of center $0$ and radius $\alpha$ of $W^{2,\infty}(\R^n,\R^n)$. This set is denoted $\m{M}$.
%
%
%
We fix an observation time $T \geq  2 (\hbox{diam}(\Omega_0))$.
%

Let $f \in L^2(\mathbb{R}^+, H^{-1}(\mathbb{R}^n))$, $y^0  \in \, H_1^0(\Omega_{\f}(0))$ and $y^1  \in \, L^2(\Omega_{\f}(0))$. We are interested in $y_{\f}$ the solution of the equations:
\begin{equation}\label{eq:WavePerturbed}
\begin{cases}
        \partial_t^2 y_{\varphi}(t,x) - \Delta y_{\varphi}(t,x) = f(t,x) & \forall t>0, x\in \Omega_{\varphi}(t)\\
        y_{\varphi}(0,x)= y^0(x)  & \forall x \in \Omega_{\f}(0)\\
        \partial_t y_{\varphi}(0,x)= y^1(x) & \forall x \in \Omega_{\f}(0)\\
        y_{\varphi}(t,x) = 0 & \forall t>0,\; x \in \partial\Omega_{\varphi}(t).
\end{cases}
\end{equation}
We denote by $y_0$ be the solution associated to the trivial perturbation $\f \equiv 0$, and $\Lambda$ the map:
\begin{equation}
\Lambda : \begin{cases}
            \m{M} \rightarrow H^1(\omega_0) \times L^2(\omega_0)\\
            \f \rightarrow  (y_{\f}(T,\cdot)\vert_{\omega}, \partial_t y_{\f}(T,\cdot)\vert_{\omega}) .
         \end{cases}
\end{equation}
The question we address is to characterize the set of traces at $t=T$ and on $\omega$ of the solutions of this problem when $\f$ is an admissible transformation, i.e. the range of the operator $\Lambda$: 
$$ \m{R}(\m{M}) = \{(y_{\f}(T,\cdot), \partial_t y_{\f}(T,\cdot))\vert_{\omega} ; \f \in \m{M}\}.$$

%


The same method as the one we used for the heat equation yields to the following equivalent variational formulation of equations~\eqref{eq:WavePerturbed} holding on $\R^+\times \Omega_0$:
\begin{multline}
\int_{\mathbb{R}^+ \times \Omega_0} \{\partial_t^2
    \bar{y}_{\varphi}|\hbox{det}(id + \nabla \varphi(t))| - \hbox{div}(B(\varphi(t)))\nabla
    \bar{y}_{\varphi}\}
    \Psi dt\,dx \\
     = \int_{\mathbb{R}^+ \times
    \Omega_0} f \circ (id + \varphi(t)) \Psi |\hbox{det}(id + \nabla \varphi(t) )| dt\,dx
\end{multline}
where 
\[ B(\varphi(t)) = \Big|\hbox{det}\big(id + \nabla \varphi(t)\, \big)\Big| \left ( \Big [\big ( \nabla(id + \varphi(t) \,) \big )^*\Big ]^{-1}\right )^* \left [\Big (\nabla\big(id + \varphi(t) \big ) \Big)^*\right ]^{-1}.\]
We observe that $B(\f)$ is symmetrical, positive because it is a perturbation of identity (this condition constrains our choice of $\alpha$ the maximal norm of $\f$ in $\m M$). We now consider $L^2(\Omega_0)$ equipped with the dot product :
\[\langle(a,b) , (c,d)\rangle_{H^1_0(\Omega_0) \times L^2(\Omega_0)} = \int_{\Omega_0}
B(\f(t)) \nabla a . \nabla c  +  \int_{\Omega_0} \,b\, d \, |\hbox{det}(id
+ \nabla \varphi (t))| dx .\]
A proof analogous to the one performed in parabolic case, proposition~\ref{theo:existenceUniquenessHeat} (based on an application of a theorem due to Lions \cite{B}) ensures existence and uniqueness of solution for the perturbed system.
%
%
%
%
%


The Lagragian shape derivative for the wave equation satisfies the equations:
\begin{equation}\label{der_lag_ondes}
     \begin{cases}
            \partial_t^2 y_{\varphi}' - \Delta y_{\varphi}' = 0\\
            y_{\varphi}'(t=0)=0\\
            \partial_t y_{\varphi}'(t=0)=0\\
            y_{\varphi}'(t)|_{\partial \Omega_0} =-\varphi(t).n \frac{\partial y_0(t)}{\partial n}.
      \end{cases}
\end{equation}



We now show that $d\Lambda(0)$ is surjective on $H^1_0(\omega) \times L^2(\omega)$, ensuring an exact controllability property of the wave equation using a local surjectivity theorem. Note that the proof we provide here differs significantly from the proof provided in the parabolic case. It is based on a controllability result due to Lions~\cite{L2} and makes use of the reversibility of the wave equation.
\begin{theorem}
 Assume that there exists $\gamma \in \partial \omega$ such that $\frac{\partial y_0}{\partial n} \geq \varepsilon >0 $ for all times $t\in [0,T]$. Then the linearized function $d\Lambda(0)$ is surjective from $\LL^2((0,T)\times \Omega_0)$ onto $H^1_0(\omega) \times \LL^2(\omega)$. 
\end{theorem}

\begin{proof}
The proof of this theorem is based on classical results on the null boundary controllability of the wave equation, proved by Lions in \cite{L1, L2}. Using Hilbert Uniqueness Method (HUM), Lions considers an open bounded subset $\Omega$ of $\RR^n$ with a smooth boundary $\partial\Omega$, a control time $T>2 \textrm{diam}(\Omega)$ and an open subset $\gamma \subset \partial\Omega$. Lions shows that for any initial condition $(u_0,u_1)\in\LL^2(\Omega)\times H^{-1}(\Omega)$, there exists a function $v\in\LL^2((0,T)\times\gamma)$ such that the solution of the equation
        $$\left \{
            \begin{array}{cccl}
                \partial_t^2 u - \Delta u &=& 0 & \textrm{ on }
                \Omega \times (0, T)\\
                u(0)&=& u_0 & \textrm{  on  } \Omega\\
                \partial_t u(0)&=& u_1 & \textrm{  on  } \Omega\\
                u(t,x) &=& v(t,x) & \textrm{  on  } \gamma \times
                (0,T)
            \end{array}
        \right .$$
\noindent is such that $u(T,x)=\partial_t u(T,x)=0$. The function $z(t,x)=u(T-t,x)$ satisfies the heat equation with zero boundary conditions and control $v(T-t,x)$ on $\gamma$, and has the property that $z(T,x)=u_0$ and $\partial_t z(T,x)=u_1$ on $\gamma$. Let us also remark that under the assumptions of the theorem, there exist several functions $\psi\in \LL^2((0,T),W^{2,\infty}(\R^n,\R^n))$ such that $-\psi(t,x).n \;\frac{\partial y_0(t,x)}{\partial n} = v(t,x)\mathbbm{1}_{x\in\gamma}$. For any of these $\psi$, we have $z=y_{\psi}'$ solution of equation~\eqref{der_lag_ondes}. Therefore, the linearized function $d\Lambda(0)$ defined by:
\begin{equation*}
d\Lambda(0) :
 \begin{cases}
   \LL^2([0,T], W^{2,\infty}(\R^n,\R^n))  &\longrightarrow  H^1_0(\omega) \times L^2(\omega) \\
  \psi &\longrightarrow (y_{\psi}'(T,\cdot)\vert_{\omega}, \partial_t y_{\psi}' (T,\cdot)\vert_{\omega})
 \end{cases}
\end{equation*}
is surjective. 
   %
\end{proof}

    The surjectivity property on the differential of $\Lambda$ directly implies the exact controllability of the original wave equation with respect to the shape of the domain, stated in the following:
	\begin{theorem}\label{theo:WaveControl}
        Let $y_0$ the solution of the unperturbed problem:
        $$ \left \{
            \begin{array}{cccl}
                \partial_t^2 y_{0} - \Delta y_{0} & =& f & \textrm{on $\Omega_{0}$}\\
                y_{0}(t=0)&=& y^0 & \textrm{ on  $\Omega_0$}\\
                \partial_t y_{0}(t=0)&=& y^1 & \textrm{ on $\Omega_0$}\\
                y_{0}|_{\partial\Omega(t)}(t) &=&0 &\forall t>0.
            \end{array}
        \right .$$
        There exists a neighborhood of $(y_0(T)\vert_{\omega}, \partial_t y_0(T)\vert_{\omega})$ in $H^1_0({\omega})
        \times L^2({\omega})$, denoted $\m{N}$ such that for all $A \in \m{N}$, there exists $\f \in \m{M}$ such that $A=\Lambda(\f)$.
    \end{theorem}
    \begin{proof}
        We proved that $\Lambda$ was differentiable and that its differential at $0$ is surjective. So the local surjectivity theorem (see e.g. Luenberger \cite{Luenberger1969} )proves theorem \ref{theo:WaveControl}.
    \end{proof}

We emphasize on the fact that the controllability property, in that case, holds both for the linearized and the original non-linear (with respect to the control $\f$) problem. Let us emphasize the fact that in the present case, the control cannot be considered constant in time in contrast with the case of the heat equation. 

\section{Control of the semi-discrete wave and heat equations}
\label{sect:semi-discrete}
In the previous section, we addressed the problem of the controllability of the heat and wave equations with respect to the shape of the domain and proved that the hyperbolic equation was locally exactly controllable with respect to the shape of the domain, the parabolic problem was not exactly controllable and its linearization was approximately controllable.

The approximate controllability property of the parabolic equation did not imply an analogous property for the nonlinear problem: though the range of the trace operator was dense, this did not imply a local inversion property because of the problem holds in an infinite-dimensional space. We shall now turn our attention to the discretized problem. This study has two main interests. First, it is relevant from a computational point of view and, second, from the mathematical point of view, since the discretized version of the Laplacian operator is finite dimensional, the density of the range of the linearized operator will imply surjectivity of the nonlinear operator.

For simplicity, this section is restricted to the analysis in two dimensions, where the open set $\Omega_0$ is a square. Using the classical methods as proposed in \cite{CZ}, it would be possible to extend these results to general domains, with an important increase of complexity in the notations, but no profound change in the mathematical arguments. Moreover, note that we are interested here in a fixed discretization of the open set $\Omega_0$. In other words, the approach does not address the convergence of this control as the stepsize of the mesh tends to zero. 

In details, we are interested in the semi-discrete heat and wave equations on a rectangle $[0,a] \times [0,b] \in \mathbb{R}^2$ discretized this set with a step $h$. The infinite continuous space system is therefore replaced by a finite-dimensional evolution problem on the discrete points of 
$$\Omega_h=\{m=(ih,jh) ; (i,j) \in \{0,\ldots,M\} \times \{0,\ldots,N\}\}.$$
Both the heat and wave equations make use of the spatial Laplacian operator, which we now define in the discretized setting. To this purpose, we introduce the following definitions:
\begin{definition} 
	Let $m=(ih,jh) \in (\mathbb{Z}h)^2$ and the discrete neighborhood of $m$ defined: 
\[B(m)=\{(kh,lh); (k,l)=(i,j),(i-1,j), (i+1,j), (i,j-1), (i,j+1)\}\]
The set of strict neighbors of $m$ is $\mathcal{B}(m)=B(m) \setminus \{m\}$.
\end{definition}

\begin{definition} The discrete interior of $\Omega_h$ is defined by 
	\[\interieur{\Omega_h} = \{m \in \Omega_h ; \mathcal{B}(m) \subset \Omega_h \}.\] 
	The discrete boundary of $\Omega_h$ is defined by $\Gamma_h = \Omega_h \setminus \interieur{\Omega_h}$ and the exterior of $\Omega_h$ as: $\interieur{F_h}=(\mathbb{Z}h)^2 \setminus \Omega_h$.
\end{definition}
%
%
These sets form a partition of $(\mathbb{Z}h)^2$. We will assume for simplicity that a single edge of the boundary is moving, for instance $\{(i,j); i=0\}$. This means that the only moving part of this set is our control. Furthermore, the free points of the boundary will move only along the normal to this boundary. 

\begin{definition}[Functional spaces]
We denote $\mathcal{F}(X)$ the set of real-valued mappings from a space $X$ and $\mathcal{F}_0(X)$ those vanishing on the boundary of $X$. 
We consider in this paper time dependent maps, taking values in $\mathcal{F}(\Omega_h)$. In particular we will use 
$C([0,T], \mathcal{F}(\Omega_h))$, the set of continuous functions $[0,T] \mapsto \mathcal{F}(\Omega_h)$ and the spaces 
$L^p([0,T],\mathcal{F}(\Omega_h))$.
\end{definition}

\begin{remark}
 The set $\mathcal{F}(\Omega_h)$ is isomorphic and identified to $\R^{MN}$. 
\end{remark}

We are now in a position to define a discretized version of the Laplacian operator, as follows.

\begin{definition}
  Let $A$ be the finite difference operator with Dirichlet boundary conditions defined by:
  \begin{equation*}
   \begin{cases}
      \mathcal{F}_0(\Omega_h) & \longrightarrow \mathcal{F}(\interieur{\Omega_h})\\
      \phi & \longrightarrow  A \phi
   \end{cases}
  \end{equation*}
   \noindent where 
    \begin{equation}\label{eq:finDifOpHeat}
        \forall m \in \interieur{\Omega_h} \textrm{,   }
        [A\phi]_m = \inv{h^2}[4 \phi(m) - \sum_{p \in B(m), p \neq
        m} \phi(p)]
    \end{equation}
\end{definition}

\subsection{Semi-discrete heat equation in a square}
\label{sect:discreteHeatEq}

We now turn our attention specifically to the case of the semi-discretized heat equation. The reference domain $\Omega_h$ being fixed, we consider a reference state $u$ as the solution of the equation
\begin{equation}\label{eq:ReferenceState}
    u \in \mathcal{F}_0(\Omega_h) \textrm{ such that }
    \left \{
        \begin{array}{ccc}
            \partial_t u + A u & = & F\\
            u(t=0) &=& u_0
        \end{array}
    \right .
\end{equation}
where the source term $F \in \mathcal{F}(\interieur{\Omega_h})$ is deduced from the source term $f$ of the continuous initial problem by a simple discretization\footnote{If $f$ is continuous, then $F$ is defined by $F_m =f(m)$ for $m \in \interieur{\Omega_h}$. If $f$ is not continuous, for instance is $f \in L^2 \textrm{ or } H^1$, then  $F_m$ will be a mean value of $f$ on a neighborhood of $m$.}.


\subsubsection{The perturbed problem}
As in the continuous case, we are interested in small perturbations of the shape of the domain $\Omega_h$. Changing the shape of $\Omega_h$ consists in moving continuously the nodes of the mesh corresponding to $x=0$. On this new subset, the finite difference Laplace operator is modified as follows:

\begin{definition}\label{def:Vj}
    We consider the set $\{V_j; j=1..N-1\}$ of vector fields $\Omega_h \mapsto \RR^2$ by: 
    \begin{equation*}
     \forall j \in \{1..N-1 \},
     \begin{cases}
      V_j(m) = (0,0) & \textrm{ if } m \neq (0, jh)\\
      V_j(m) = (1,0) & \textrm{ if } m = (0, jh).
     \end{cases}
    \end{equation*}
    Let $W_h$ be the vector space spanned by the family $(V_j)_{j \in \{1...N-1\}}$. The perturbations we consider in this problem are in the set:
    \begin{multline}
     \mathcal{W}_h=\{ \sum_{j=1}^{N-1} h \, \lambda_j(t) V_j ; t\rightarrow \lambda_j(t) \in L^{\infty}(\RR^+;W_h)\cap C(\RR^+;W_h) \\
     \textrm{ and such that  }\; \sup_{j=1..N-1} \|\lambda_j\|_{\infty} < 1/2
    \}
    \end{multline}
\end{definition}

\begin{remark}
Note that the perturbation has the same magnitude as $h$. This does not allow us even in the better cases to have the continuous case as a limit case since the perturbation tends to the trivial condition as the mesh becomes finer.
\end{remark}

\begin{definition}
    Let us define $\Gamma^1$ the first layer of interior nodes:
    \[\Gamma^1=\{(1,j)\textrm{   ;   } j\in\{1,..N-1\}\}.\]
\end{definition}

The perturbed heat operator on this new discretized set is defined as:
\begin{definition}\label{def:PerturbedOp}
    Let $\varphi(t)= \sum_{j=1}^{N-1} \lambda_j(t) h V_j \in \mathcal{W}_h$. The operator $A(\varphi) : \m{F}_0(\Omega_h^{\varphi}) \mapsto \m{F}(\interieur{\Omega_h})$ is defined as:
    \begin{equation}
    \begin{cases}
     \inv{h^2}[4 \phi(m) - \sum\limits_{p \in B(m),\; p \neq m} \phi(p)] & \forall m \in \interieur{\Omega_h} \setminus \Gamma^1 \\
     \inv{h^2}[2(1+\inv{1+\lambda_j(t)})\phi_{(1,j)}-\frac{2}{2+\lambda_j(t)}\phi_{(2,j)}- \phi_{(1,j+1)}-\phi_{(1,j-1)}] & \textrm{ for } m=(1,j) \in \Gamma^1
    \end{cases}
    \end{equation}
\end{definition}

\begin{proposition}\label{prop:Abounded}
    The operator $A(\varphi)$ is bounded for all $\varphi \in \m{W}_h$.
\end{proposition}

\begin{proof}
	It is easy to show that:
    \begin{itemize}
        \item For $m \in \interieur{\Omega_h} \setminus \Gamma^1$, we have $[A(\varphi)\phi]_m = [A\phi]_m=\inv{h^2}[4 \phi(m) - \sum_{p \in B(m), p \neq
            m} \phi(p)]$, and hence:
            \begin{equation*}
                |[A(\varphi)\phi]_m| \leq \inv{h^2} \left [ 4 \|\phi\|_{\infty} + \sum_{p \in B(m), p \neq m} \|\phi\|_{\infty} \right ] \leq \frac{8}{h^2} \|\phi\|_{\infty}.
            \end{equation*}
        \item For $m = (1,j) \in \Gamma^1$, we have 
            \[[A(\varphi)\phi]_m=\inv{h^2}(2(1+\inv{1+\lambda_j(t)})\phi_{(1,j)}-\frac{2}{2+\lambda_j(t)}\phi_{(2,j)}
            - \phi_{(1,j+1)}-\phi_{(1,j-1)}),\] and hence we have:
            \begin{equation*}
                |[A(\varphi)\phi]_m| \leq \inv{h^2}\left (2(1+2)\|\phi\|_{\infty}+\frac{2}{2-1/2}\|\phi\|_{\infty}
                + 2 \|\phi\|_{\infty} \right ) \leq \frac{28}{3 h^2} \|\phi\|_{\infty}.
            \end{equation*}
    \end{itemize}
    Since we are in $\m{F}(\interieur{\Omega_h})$ which is a finite dimension vector space, all the norms are equivalent so the operator $A(\varphi)$ is bounded on $\m{W}_h$.
\end{proof}



\begin{definition}
    The perturbed state $\overrightarrow{u_{\varphi}}(x,t) \in \m{F}(\interieur{\Omega_h})$ is the unique solution in  $C([0,T],
    \m{F}(\interieur{\Omega_h}))$ of the semi-discrete problem:
    \begin{equation}\label{eq:def-etat-pertube}
        \begin{cases}
         \partial_t \vect{u_{\varphi}} + A(\varphi)\vect{u_{\varphi}}&= F\\
         \vect{u_{\varphi}}(t=0) &= \vect{u_0}.
        \end{cases}
    \end{equation}
    This solution exists and is unique, and defined for all time $t>0$.
\end{definition}

\begin{proof}
	    Equation~\eqref{eq:def-etat-pertube} is a linear ordinary differential equation with $t$-measurable vector field, hence classical theory (Cauchy-Lipschitz theorem) implies local existence and uniqueness of the perturbed state, which will be continuous in its definition domain. Non-explosion property is a classical application of Gronwall's lemma based on the boundedness of the perturbed operator and of the source term $F$ (proposition \ref{prop:Abounded}).
	%
\end{proof}

\subsubsection{Controllability of the semi-discrete heat equation}\label{sec:ControlSemiDiscreteHeat}

As in the continuous case, we consider the map
\begin{equation}\label{eq:lambdah}
    \Lambda_h : \begin{cases}
            \m{W}_h &\longrightarrow \m{F}(\interieur{\Omega_h})\\
            \varphi &\longrightarrow \vect{u_{\varphi}}(T)
        \end{cases}
\end{equation}

\noindent where $\vect{u_{\varphi}}(T)$ is solution of equation \eqref{eq:def-etat-pertube}. Let $Z_d = u(T)$, where $u$ is the reference state \eqref{eq:ReferenceState}. The problem we address is to find a neighborhood $\m{V} (0) \in \m{V}$ of the reference domain and $\m{V}(Z_d) \in \m{F}(\interieur{\Omega_h})$ of the trace at $t=T$ of the reference solution such that $\m{V}(Z_d) \subset
\Lambda_h(\m{V}(0))$.\\

In the finite-dimension spaces where the problem is now set, we use the local inversion theorem to demonstrate this property. 
First of all we will prove that $\Lambda_h$ is differentiable in the neighborhood of the origin, and that $d\Lambda_h(0)$ is surjective. Then we will use the adjoint state technique to prove that the surjectivity of $d\Lambda_h(0)$ is equivalent to a pool of conditions the semi-discrete adjoint should state satisfies. Finally, we will prove the controllability property proving those conditions on the adjoint state, which happen to be a property of discrete unique continuation.

{\bf (i). Differentiability}\\

Let us denote $tr$ the trace operator:
\[tr:\begin{cases}C([0,T],\m{F}(\Omega_h)) &\mapsto \R \\ \phi &\mapsto \phi(T) \end{cases} \] 
The map $\Lambda_h$ is the composition of the trace operator and the map $U : \varphi \longrightarrow
\vect{u_{\varphi}}$.

The trace function is linear and continuous. So we only need to prove that $U$ is Fr\'echet-differentiable in $0$.

\begin{proposition}\label{prop:existencegateaux}
    $\forall \f \in \m{W}_h \textrm{ ,  } \forall \psi \in W_h$, the map $ \f \mapsto \vect{u_{\f}}$ is 
    differentiable in $\f$ in the direction of $\psi$, and the G\^ateaux differential $\langle D_G \Lambda_h(\varphi), \psi\rangle$, denoted $\vect{v_{\varphi}}(\psi)$ is solution of the differential equation: 
    \begin{equation}\label{eq:DerGateaux}
      {\begin{cases}
            \partial_t \vect{v_{\varphi}}(\psi) + A(\varphi)\vect{v_{\varphi}}(\psi) &= -\langle\dot{A_{\varphi}}, \psi\rangle   \vect{u_{\varphi}}\\
            \vect{v_{\varphi}}(\psi) (t=0) &= \vect{0}.
      \end{cases}}
    \end{equation}
\end{proposition}

\begin{proof}
First let $\displaystyle{W_{\lambda}=\frac{u_{\varphi + \lambda \psi}-u_{\varphi}}{\lambda}}$. The function $D_G \Lambda_h$ is the limit, when it exists, of $W_{\lambda}$ when $\lambda$ tends to $0$.

\begin{enumerate}
    \item {\it Necessary condition}: We assume that this limit exists, denote it $v$, and we compute the equation this limit satisfies.
    The function $u_{\varphi + \lambda \psi}$ satisfies the equations:
    \begin{equation*}
     \begin{cases}
      \partial_t \vect{u_{\varphi + \lambda \psi}} + A(\varphi + \lambda \psi) \vect{u_{\varphi + \lambda \psi}} &= F(t)\\
      \vect{u_{\varphi + \lambda \psi}} (t=0) &= \vect{u_0}.
     \end{cases}
    \end{equation*}

    We differentiate this equation with respect to $\lambda$ at $\lambda=0$, and we obtain the equation $v$,
    \begin{equation*}
     \begin{cases}
      \partial_t v + \langle  \dot{A}({\varphi}), \psi\rangle  \vect{u_{\varphi}} + A(\varphi) \vect{v} &= 0\\
      \vect{v}(t=0) &= 0.
     \end{cases}
    \end{equation*}

    So we deduce that if the differential of $W_{\lambda}$ exists, then it is the solution $v_{\varphi}(\psi)$ of the ordinary differential equation \eqref{eq:DerGateaux}:\\
    \begin{equation*}
    \begin{cases}
            \partial_t \vect{v_{\varphi}}(\psi) + A(\varphi)\vect{v_{\varphi}}(\psi) &= -\langle  \dot{A}({\varphi}), \psi\rangle   \vect{u_{\varphi}}\\
            \vect{v_{\varphi}}(\psi) (t=0) &= \vect{0}.
        \end{cases}
    \end{equation*}

    \item {\it Sufficient condition}: We show that the solution of the ordinary differential equation \eqref{eq:DerGateaux} is the limit of $W_{\lambda}$.
    Indeed, the map $(\lambda_j)_{j= 1 ... N-1} \mapsto A(\phi)$ is $C^{\infty}$, so the differential $\langle  \dot{A_{\varphi}},\psi\rangle  $ is defined, and furthermore $\langle  \dot{A_{\varphi}},\psi\rangle  \vect{u_{\varphi}}$ is $L^2(0,T)$. So Cauchy-Lipschitz theorem ensures existence and uniqueness of solution, and we can prove that it is defined for all time using Gronwall's lemma. Let us now show that the function $\vect{v}(t)$ defined is indeed the limit of $W_{\lambda}$ when $\lambda \to 0$. Let $W_{\lambda} - v := \varepsilon_{\lambda} $. We have: 
    \begin{align*}
    \varepsilon_{\lambda} &= \partial_t \varepsilon_{\lambda} + A(\varphi + \lambda \psi) W_{\lambda}- A(\varphi) v \\
    & = \left [\frac{A(\varphi + \lambda \psi)-A(\varphi)}{\lambda}-\langle  \dot{A}({\varphi}),\psi\rangle  \right ]u_{\varphi}.
    \end{align*}
    Since $A \in C^2$ in $\varphi$, the right hand of the equality is $\m{O}(\lambda)$, as well as $(A(\varphi + \lambda \psi)-A(\varphi))v$ which comes from the left hand of the equality. So eventually, $\varepsilon_{\lambda}$ satisfies the equation:\\
    \begin{equation*}
     {
     \begin{cases}
      \partial_t \varepsilon_{\lambda} + A(\varphi)\varepsilon_{\lambda} &= o(\lambda)\\
      \varepsilon_{\lambda}(t=0) &= \vect{0}
     \end{cases}
     }
    \end{equation*}
	It is hence clear that
    $\|\partial_t \varepsilon_{\lambda}\| \leq \alpha \|\varepsilon_{\lambda}\| + M\lambda$
    and hence $\|\varepsilon_{\lambda}(t)\|\leq \lambda \frac{M}{\alpha}e^{\alpha T}$. Therefore, $\varepsilon_{\lambda}$ converges uniformly to $0$ when $\lambda \to 0$ ensuring that the limit of $W_{\lambda}$ exists and is indeed the function $v$ solution of \eqref{eq:DerGateaux}.
\end{enumerate}
\end{proof}

%
%
%
%

\begin{proposition}\label{deriveecontinue}
The differential of $u_{\phi}$ with respect to $\phi$, $$\varphi
\longrightarrow D_G \Lambda_h(\varphi)$$ is continuous.
\end{proposition}

\begin{proof}
$D_G \Lambda_h(\varphi) : \psi \longrightarrow \vect{v}$ where $\vect{v}$ is solution of:
$$ \left \{
    \begin{array}{ccc}
        \partial_t \vect{v} + A_{\varphi} \vect{v} &=&
        -\langle \dot{A_{\varphi}},\psi\rangle \vect{u_{\varphi}}\\
        \vect{v}(t=0)&=&0.
    \end{array}
    \right .
    $$
Hence $\vect{v}$ is solution of $\partial_t \vect{v} =
F_{\varphi, \psi}(t,\vect{v})$ with $F_{\varphi,
\psi}(t,\vect{x}) = -A_{\varphi} \vect{x}
-\langle  \dot{A_{\varphi}},\psi\rangle   \vect{y_{\varphi}}.$ We note that the function $$\varphi \rightarrow F_{\varphi,
\psi}$$ is continuous, since $\varphi \rightarrow A_{\varphi}$ is continuous. Cauchy-Lipschitz' theorem with parameters gives us that $\varphi \rightarrow y_{\varphi}$ is also continuous. 

Moreover, the differential $\langle  \dot{A_{\varphi}},\psi\rangle  $ is continuous in $\f$ since $A_{\varphi}$ has a rational variation in $\varphi$, and by definition of $\varphi$, these rational fractions have no singular point on $\phi$, implying that  the dependence in $\varphi$ remains continuous.

So Cauchy-Lipschitz' theorem with parameters gives us the continuity of $y_{\varphi}'$ in $\varphi$.
\end{proof}

\begin{theorem}
$\Lambda_h$ is Fr\'echet-differentiable at $0$.
\end{theorem}

\begin{proof}
We already proved that:
\begin{itemize}
    \item the G\^ateaux differentials in all directions of $ \m{V}$ exist (Prop.\ref{prop:existencegateaux})
    \item these differentials are continuous (Prop. \ref{deriveecontinue})
\end{itemize}
Using the property that a function G\^ateaux differentiable in all directions and the differential of which being continuous is differentiable, we conclude on the Fr\'echet differentiability of  $\Lambda_h$ at $0$.
\end{proof}

\noindent{\bf (ii). Adjoint state technique}\\
We recall that $\varphi \in L^2(0,T;W_h)$. We have: 
$$\varphi(t)=\sum_{j=0}^{N-1} h \lambda_j(t) V_j.$$ 

\noindent $\Lambda_h$, defined in \eqref{eq:lambdah} is the composed application of the trace function at $t=T$ 
%
with the map $U$ defined by: 
$$U :
\left \{
    \begin{array}{lcl}
        C([0,T]; W_h) &\longrightarrow &
        C([0,T];\m{F}(\interieur{\Omega_h}))\\
        \varphi &\longrightarrow & u_{\varphi} 
    \end{array}
\right .
$$

\noindent where $u_{\varphi}$ is solution of the equations:
\[
 \left \{
  \begin{array}{ccc}
    \partial_t u_{\varphi} + A(\varphi) u_{\varphi}
      &=& F\\
      u_{\varphi}(t=0)&=&0.
\end{array}
\right .
\]
Furthermore, recall that $$d\Lambda_h(0) : \left \{
    \begin{array}{ccc}
        C([0,T], W_h) &\longrightarrow
        &\m{F}(\interieur{\Omega_h})\\
        \varphi &\longrightarrow  & y_{\varphi}(T)
    \end{array}
    \right .
    $$
    and  $y_{\varphi}$ is solution of
    $$
    \left \{
    \begin{array}{ccc}
        \partial_t y_{\varphi} + A(0) y_{\varphi} &=& -\langle \dot{A}(0),
        \varphi \rangle u_0\\
        y_{\varphi} (t=0) &=&0.
    \end{array}
    \right .
    $$
Eventually, we denote $Y : \varphi \longrightarrow y_{\varphi}$ (so we have $d\Lambda_h (0) = tr|_{t=T} \circ Y$) and remark that the adjoint of the trace map is given, for any $c \in \m{F}(\interieur{\Omega_h})$, by $tr^*c = c \delta_{t=T}$. We now prove that the map $d\Lambda_h (0)$ is surjective. To this purpose, we use the adjoint state technique. 

The surjectivity of $d\Lambda_h (0)$ is equivalent to the fact that:


\[\{c \in \m{F}(\interieur{\Omega_h});
    \forall \varphi \in \m{W}_h \langle d\Lambda_h (0) \varphi, c\rangle  =0\}= \{0\}.\]

Any $c \in \m{F}(\interieur{\Omega_h})$ such that $\forall \varphi \in \m{W}_h \langle d\Lambda_h (0) \varphi, c\rangle  =0$ is such that $\langle tr|_{t=T}(Y(\varphi)),c\rangle  =0$ for any $\varphi \in \m{W}_h$, which is equivalent to the property:
    
    \begin{eqnarray}\label{condition_surjective}
    \forall \varphi \in
    \m{W}_h \langle Y(\varphi),tr|_{t=T}^*(c)\rangle =0
    \end{eqnarray}

\begin{definition}\label{def_etat_adjoint}
    The adjoint state associated to $y_{\varphi}$ and $tr|_{t=T}^*(c)$ is the unique solution $X$ of the equations:
    \begin{equation}\label{defadj}
    \begin{cases}
            -\partial_t X + A X &= tr|_{t=T}^*(c)\\
            X(t=0)&=0
    \end{cases}
    \end{equation}
\end{definition}

\begin{remark} Using the fact that $A$ is self-adjoint, we clearly have:
\begin{eqnarray}
    \nonumber \langle \partial_t y_{\varphi} + A y_{\varphi}, T\rangle   &=& \langle \partial_t y_{\varphi}, T \rangle  + \langle A y_{\varphi},T \rangle  \\
    \nonumber &=& -\langle y_{\varphi},\partial_t  T \rangle + \langle y_{\varphi},A^*T \rangle \\
    \nonumber &=& \langle y_{\varphi}, -\partial_t  T + AT \rangle
\end{eqnarray}
\end{remark}

Using this definition we replace in \eqref{condition_surjective} $tr|_{t=T}^*(c)$ by its expression in function of the adjoint state $X$ defined in \eqref{defadj} and obtain the set of equivalent statements:
$$
    \begin{array}{lll}
        &\forall \varphi \in \m{W}_h &
        \langle Y(\varphi),tr|_{t=T}^*(c)\rangle =0\\
        \Leftrightarrow &\forall \varphi \in \m{W}_h &
        \langle Y(\varphi),-\partial_t X + A X \rangle =0\\
        \Leftrightarrow & \forall \varphi \in \m{W}_h &
        \langle \partial_t y_{\varphi} + A y_{\varphi},X\rangle =0\\
        \Leftrightarrow & \forall \varphi \in \m{W}_h &
    \langle  A'_{\varphi}y_0,X\rangle   =0.
\end{array}
$$

This proves the following:

\begin{theorem}\label{surjection_theoreme}
    The differential $d\Lambda_h(0)$ of  $\Lambda_h$ at $\varphi = 0$ is surjective if and only if we have the following uniqueness property:  If $c\in \m{F}(\interieur{\Omega_h})$ is such that 
    \begin{equation}\label{condi_surj_theo}
    \langle  X, A'_{\varphi}y_0 \rangle_{L^2(0,T;\m{F}(\interieur{\Omega_h}))} =
    0, \textrm{   } \forall \varphi \in \m{W}_h
    \end{equation}
    \noindent where X is solution of :
    \[\left \{
            \begin{array}{ccc}
                -\partial_t X + A X &=& tr|_{t=T}^*(c)\\
                X(t=0)&=&0,
            \end{array}
        \right .
    \]
    then necessarily $c=0.$
\end{theorem}

\noindent{\bf (iii). Calculation of the differential of $A$ at $0$}

\begin{proposition}\label{aprim}
    Let  $j \in \{1...N-1\}$ and $\phi \in \fom$. For all $\mu \in C([0,T],\R)$, we denote $\langle A'_{0},\mu(t)V_j\rangle $ the differential of $A$ at $0$ in the direction $\mu(t)V_j$. We have:
    \begin{equation}
    [\langle  A'_{0}, \mu(t)V_j\rangle \phi]_m  =
     \begin{cases}
       0 & \forall m \in \interieur{\Omega_h} \setminus \Gamma^1\\
       \frac{\mu(t)}{h^2} \left (\inv{2} \phi_{(2,j)} - 2 \phi_{(1,j)}\right) & \textrm{ if }\; m=(1,j)
     \end{cases}
    \end{equation}
\end{proposition}

\begin{proof}
 For all $j \in \{1...N-1\}$, the point $(0,j)$ of the boundary has a unique neighbor in $\Om$, which is $(1,j)$. From the definition \ref{def:PerturbedOp} of $A(\varphi)$, for all $\mu : t  \rightarrow \mu(t)$ $C([0,T]; \RR)$, we have :
 $$\forall m \in \interieur{\Omega_h} \textrm{ ,  }
        m \neq (1,j) \textrm{ ,  } [A_{\mu(t)V_j}\phi]_m = [A\phi]_m$$
 and we have also 
 $$ [A(\mu(t)V_j)\phi]_{(1,j)}=\inv{h^2}\left ( 2(1+\inv{1+\mu(t)})\phi_{(1,j)}-\frac{2}{2+\mu(t)}\phi_{(2,j)}
        - \phi_{(1,j+1)}-\phi_{(1,j-1)}\right).$$
 Hence we get
 \begin{eqnarray}
    \nonumber [A(\mu(t)V_j)\phi]_{(1,j)}-[A\phi]_{(1,j)}&=&
        \inv{h^2}[2(1+\inv{1+\mu(t)}-2)\phi_{(1,j)}
        -(\frac{2}{2+\mu(t)}-1)\phi_{(2,j)}].\\
    \nonumber  &=& \inv{h^2}[2(1+1-\mu(t) + o(\mu)-2)\phi_{(1,j)}
        -(1-\inv{2}\mu(t)-1+o(\mu))\phi_{(2,j)}].
 \end{eqnarray}
\end{proof}

\noindent{\bf (iv). The condition $X|_{\Gamma^1}=0$} 

\begin{definition}\label{NDD}
    The function $Y\in C([0,T], \fom)$ satisfies the discrete non-degeneracy condition if and only if:
    $$\forall t>0, \forall j \in \{1,...,N-1\},  \inv{2}Y_{(2,j)}-2Y_{(1,j)} \neq 0.$$
\end{definition}

\begin{remark}
    This condition can be seen as an finite difference approximation of the condition $\frac{\partial y}{\partial n}\neq 0$. Indeed, let $ y \in \m{C}^2(\m{V}(x))$ where $\m{V}(x)$ is a neighborhood of $x$. Assume that $y$ satisfies $y(x)=0$.
    Then performing a Taylor expansion, we get: $$y'(x)=\inv{h} [2 y(x+h)- \inv{2} y(x+2h)] + o(h).$$
    Note also that in the continuous case, we only need to assume that $\frac{\partial y}{\partial n}$ does not vanish on an open set of the boundary, and not all along the boundary. Here we need to assume the non-degeneracy condition all along the boundary of $\Omega_h$ to prove the discrete unique continuation.
\end{remark}

\begin{proposition}\label{layer1}
    Under the discrete non-degeneracy condition (definition~\ref{NDD}) on the reference state, we have:
    $$X|_{\Gamma^1} \equiv 0.$$
\end{proposition}

\begin{proof}
The relation \eqref{condi_surj_theo} gives us  : 
$$\langle  X, A'_0(\mu(t)V_j)y_0\rangle_{L^2(0,T);\m{F}(\interieur{\Omega_h})} = 0, \textrm{   } \forall j \in \{1,..,N-1\}, \forall \mu(t) C_0^{\infty}(\RR^+, \RR).$$
Proposition \ref{aprim} gives us therefore that $\forall j \in
\{1,..,N-1\}$,
 \begin{itemize}
        \item $\forall m \in \interieur{\Omega_h} \; , \;
        m \neq (1,j) \; , \; [A'_0(\mu(t)V_j)\phi]_m = 0,$
        \item $[A'_0(\mu(t)V_j)\phi]_{(1,j)} = \inv{h^2}[\inv{2}
        \phi_{(2,j)} - 2 \phi_{(1,j)}]\mu(t)$
 \end{itemize}

so we eventually have : $$\langle  X, A'_0(\mu(t)V_j)y_0\rangle_{L^2(0,T);\m{F}(\interieur{\Omega_h})} = 0,$$
    Thus
    $$\int_0^T \sum_{m\in \Om} \mu(t) [A'_0(\mu(t)V_j)y_0]_m(t)
        X_m(t) dt = \int_0^T \mu(t)\inv{h^2}[\inv{2}
        y_0(2,j) - 2 y_0(1,j)] X_{(1,j)}(t) dt = 0.$$
        Therefore, for all $t>0$, $\inv{h^2}[\inv{2} y_0(2,j) - 2 y_0(1,j)](t) X_{(1,j)}(t)=0$
        
        Since $\inv{2} y_0(2,j) - 2 y_0(1,j)$ never vanishes, we have:
        $$X_{(1,j)}(t) = 0 \;\; \forall t>0 \textrm{ and  } \forall j \in
        \{1,...,N-1\}.$$
\end{proof}

\noindent{\bf (v). Unique discrete continuation} \\

The aim of this section is to prove that the uniqueness condition appearing in theorem~\ref{condi_surj_theo} (equation \eqref{condi_surj_theo}) is valid. This uniqueness condition is proved using the fact that under the discrete non-degeneracy condition on the reference state $y_0$, $X|_{\Gamma^1} \equiv 0$ (proposition \ref{layer1}). We now show that this implies that the adjoint state $X$ given by equations~\eqref{def_etat_adjoint} identically null, so $tr_{t=T}^* (c) =0$ and $c=0$. 

The method we use is based on the study of the propagation of the zeros of $X$ on $\Omega_h$ from its boundary, analogous to the approach developed in the continuous case using Holmgren's theorem. The main difference is that the propagation of zeros in the continuous case is a global property, whereas it is a local property in the discrete case \footnote{This is why we need to assume the non-degeneracy condition all along the the boundary of $\Omega_h$.}.

\begin{theorem}\label{cegal0}
    The unique solution of the equations 
    \begin{equation}\label{adjoint}
     \left \{
        \begin{array}{ccc}
            -\partial_t X + A X &=& tr|_{t=T}^*(c)\\
            X(t=0)&=&0\\
            X|_{\Gamma^0} &=& 0\\
            X|_{\Gamma^1} &=& 0
        \end{array}
    \right .
    \end{equation}
    is $X \equiv 0$, and so $c = 0$.
\end{theorem}

\begin{proof}
    \begin{enumerate}
        \item First we are interested in the equation $-\partial_t X + A X = tr|_{t=T}^*(c)$. On the sets $[0,T[$ and $]T,\infty[$, the equation simply reads $-\partial_t X + A X = 0$ so we have existence, uniqueness and continuity of the solution $X$ in these domains. The right hand term can be interpreted as an imposed jump condition at time $t=T$. Indeed, let us write the variational formulation of the problem\eqref{adjoint}: \\Let $v \in C_0^{\infty}(\RR^+,\fom)$. The variational formulation reads:
        $$\int_0^{\infty} \sum_{m\in \Om} -\partial_t X _m v_m + (AX)_m.v_m dt= \sum_{m\in \Om}v_m(T) c_m$$ i.e :
        $$\int_0^{\infty} \sum_{m\in \Om} \partial_t X _m v_m + (AX)_m.v_m dt= \sum_{m\in \Om}v_m(T) c_m$$
        The imposed jump at $t=T$ reads:
            \begin{eqnarray}
                \nonumber -\partial_t X + A X &=& 0\\
                \nonumber X(t=0)&=&0\\
                \nonumber [|X|] (T) &=& c\\
                X|_{\Gamma^0} &=&0
            \end{eqnarray}
        where we denoted $[|X|] (T)=X(T^+)-X(T^-)$ the jump of $X$ at $t=T$. The variational formulation of the problem reads: $\forall v \in C^{\infty}_0(\RR^+, \fom)$
        \begin{eqnarray}
            \nonumber 0 &=& \int_0^{\infty} \sum_{M\in \Om} -\partial_t X _m v_m + AX_m.v_m dt\\
            \nonumber &=& \int_0^{T} \sum_{M\in \Om} -\partial_t X _m v_m + AX_m.v_m dt + \int_T^{\infty} \sum_{M\in \Om} -\partial_t X _m v_m + AX_m.v_m dt\\
            \nonumber &=& \int_0^{T} \sum_{M\in \Om}  X _m \partial_t v_m + AX_m.v_m dt + \int_T^{\infty} \sum_{M\in \Om} X _m \partial_t v_m + AX_m.v_m dt \\
            \nonumber & &- \sum_{M\in \Om} (X_m(T^+)-X_m(T^-))v_m(T) + X_m(0)v_m(0)\\
            \nonumber &=& \int_0^{\infty} \sum_{M\in \Om} X _m \partial_t v_m + AX_m.v_m dt - \sum_{M\in \Om} [|X_m|](T) v_m(T).
        \end{eqnarray}
        We have the same variational formulations, so the solution are identical.
    \item Calculation of the solution $X$: we have $\forall j \in \{0...N\}, X_{0,j}=X_{1,j}=0$. We reason by induction on $k$. Assume that on the column $k-1$ ($\{(i,j);\; i=k-1\}$) and the column $k$ ($\{(i,j); i=k\}$) we had $X=0$. In this case, $X$ also vanishes on the column $k+1$.
    
    Indeed, let $j \in \{0,.., N\}$
    \begin{enumerate}
      \item If $j=0$ or  $j=N$ then we have indeed $X_{(k+1,j)}=0$ because $X$ vanishes on $\Gamma_h$, the boundary of $\Omega_h$.
      \item If $j \in \{1,.., N-1\}$. Let us write the equation satisfied by $X_{(k,j)}$:\\
      \begin{equation}\label{annulation}
         \begin{cases}
          \partial_t X_{(k,j)} + (AX)_{(k,j)}= 0\\
          \partial_t X_{(k,j)} + \inv{h^2}\left [4 X_{(k,j)}- X_{(k+1,j)} - X_{(k-1,j)} - X_{(k,j+1)} - X_{(k,j-1)} \right ] = 0.

         \end{cases}
       \end{equation}
       Since we assumed that: $X_{(k,i)}(t) = X_{(k-1,i)}(t)\equiv 0$. We reinject this condition in \eqref{annulation} and we get:
                $$-X_{(k+1,j)} (t) = 0 \;\;\; \forall t.$$
    \end{enumerate}
    So we prove that if the solution vanishes on two consecutive columns, then the solution vanishes on all the other columns. The hypothesis being that the solution $X$ vanishes on the column $i=0$ and $i=1$, we indeed proved that $X$ is identically null. 
    \item We conclude that the solution $X$ of this problem is time continuous, (it is constant equal to 0) and that the jump of the solution at $t=T$ is null, so $c=0$.
    \end{enumerate}
\end{proof}

\noindent{\bf (vi). Discrete controllability result} \\

\begin{theorem}
    Assume that the reference state $y_0$ defined by \eqref{eq:def-etat-pertube} satisfies the non-degeneracy discrete condition \eqref{NDD}. 
    Let $y_{\varphi}$ be the solution of the perturbed state \eqref{eq:def-etat-pertube} and $Z_d = y_0(T)$. 
    
    There exist neighborhoods $\m{V}(0) \subset C([0,T];W_h)$ and $\m{V}(Z_d) \subset \fom$ such that for all $Z \in \m{V}(Z_d)$ there exists $\varphi \in \m{V}(0)$ such that $y_{\varphi}(T)=Z$.
\end{theorem}

\begin{proof}
    This is a consequence of the local surjectivity property of the map $\Lambda_h$ (defined in \eqref{eq:lambdah}). 
    %
    %
    Proposition \ref{layer1} shows that when the reference state satisfies the discrete non-degeneracy condition, then $X|_{\Gamma^1} = 0$. In theorem \ref{cegal0} we proved that this second relation implies that $X=0$ and that $c=0$. This readily implies that $d\Lambda_h(0)$ is surjective using theorem~\ref{condi_surj_theo}, which completes the proof.
\end{proof}

Therefore, we have proved that the semi-discrete heat equation was locally exactly controllable, which was not the case of the continuous-space equation. We now turn our attention to the case of the wave equation. We realize again that a control independent of the time can be found in that case.


\subsection{Semi-discrete wave equation in a square}
\label{sect:discreteWaveEq}

The controllability of the waves equations is demonstrated in an analogous manner.
%
%
%
%
Similarly to the parabolic case,  the wave equation unperturbed state $u$ is solution of the ordinary differential equation: 
\begin{equation}\label{defv0_ondes}
  \begin{cases}
    u &\in \mathcal{F}_0(\Omega_h) \\
    \partial_t^2 u_0 + A u_0 &= F\\
    u_0(t=0)&=u_0\\
    \partial_t u_0 (t=0)&=u_1
  \end{cases}
\end{equation}
Where the discrete Laplace operator $A$ and the function $F$ are defined as in section~\ref{sec:ControlSemiDiscreteHeat}. 


\subsubsection{The perturbed state}
As in the continuous state, we are interested in small perturbations of the shape of the domain $\Omega_h$. We are quite free in the choice of the admissible transformations, and only look for sufficient conditions for the exact controllability. The first assumption we make on the perturbation is that the shape of the domain will be modified only moving nodes of the mesh using $C^1$ transformations in time. Moreover, the only moving nodes are located on the line $x=0$, and will move along the normal to the boundary. Because of the finite propagation speed of information in the wave equation, the problem will be well posed if the boundary does not moves faster than the information, i.e. the differential of the deformation should not have a module greater than the information propagation speed, in our case $1$.

On those perturbed open sets, the operator approximating the Dirichlet Laplacian is identical to the one defined for the heat equation.
We recall that $W_h$ is the real vector space spanned by $(V_j)_{j \in \{1...N-1\}}$ defined in~\ref{def:Vj}. The admissible transformations we consider belong to the space:
    \begin{multline*}
      \mathcal{V}=\Big\{ \sum_{j=1}^{N-1} h \lambda_j(t) V_j ;\;\; t
      \rightarrow \lambda_j(t) \in W^{1,\infty}(\RR^+;Wh)\cap C^1(\RR^+;Wh)\\
      \textrm{ and such that   } \sup_{j=1..N-1} \|\lambda_j\|_{\infty} <
      1/2, \;\; \|\partial_t \lambda_j\|_{\infty} <  1\Big \}.
    \end{multline*}

We denote by $M(\f)$ denote the $2n \times 2n$ matrix:
    $$\left (
        \begin{array}{cc}
            0 & -id\\
            A(\f) & 0
        \end{array}
    \right )
    $$

    The perturbed state, denoted $\overrightarrow{u_{\varphi}}(x,t) \in \m{F}(\interieur{\Omega_h})$, is the unique solution in $C([0,T],\m{F}(\interieur{\Omega_h}))$ of the semi-discrete problem:
    \begin{equation}\label{def-etat-pertube_ondes}
        \begin{cases}
          \partial_t^2 \vect{u_{\varphi}} + A(\varphi) \vect{u_{\varphi}}&= F\\
          \vect{u_{\varphi}}(t=0) &= \vect{u_0}\\
          \partial_t \vect{u_{\varphi}}(t=0) &= \vect{u_1}
        \end{cases}
    \end{equation}

A direct application of standard theory of ordinary differential equations ensures that:
\begin{proposition}
    We define $U_{\f}:=\left ( \begin{array}{c} u_{\varphi}\\ \partial_t u_{\varphi}  \end{array} \right )$. For all $\f \in \mathcal{V}$, $U(\f)$ is well defined and bounded in $C([0,T],\m{F}(\interieur{\w_h}))^2$
\end{proposition}

%
%
%
%

\subsubsection{Controllability of the semi-discrete wave equation}
We show a surjectivity property of the map:
\begin{equation}\label{lambdah_ondes}
    \Lambda_h : 
      \begin{cases}
            \m{V} &\mapsto \m{F}(\interieur{\Omega_h})\\
            \varphi &\mapsto U_{\varphi}(T)
      \end{cases}
\end{equation}
\noindent where $U_{\varphi}(T)$ is solution of the equation \eqref{def-etat-pertube_ondes}. Let $Z_d = u_0(T)$, where $u_0$ is the reference state defined in \eqref{defv0_ondes}. 
%

\noindent {\bf (i). Differentiability }\\
 
We denote here again $tr$ the trace function at $t=T$. The map $\Lambda_h$ is the composition of the function $S : \varphi \longrightarrow U_{\varphi}$ and the trace function. The Fr\'echet-differentiability of $\Lambda_h$ is equivalent to the Fr\'echet-differentiability of $S$ at 0. 

This differentiability is an immediate consequence of the differentiability of $A(\f)$. Indeed, by
Cauchy-Lipschitz' theorem with parameters, if $A$ is Fr\'echet-differentiable in $\f$, then the matrix 
$$M(\f)=
\left (
        \begin{array}{cc}
            0 & -id\\
            A(\f) & 0
        \end{array}
    \right )
$$
is differentiable in $\f$ which gives the Fr\'echet-differentiability of  $U_{\f}$ in $\f$.

\begin{proposition}
    $\phi \rightarrow U_{\phi}$ is differentiable at $0$ in the direction $\psi$, and the Fr\'echet-differential  $\langle  d\Lambda_h(0), \psi\rangle$, denoted $Y_{\psi}$, is solution of the differential equation: 
    $$
    \left \{
        \begin{array}{lll}
            \partial_t Y_{\psi} + M(0)
            Y_{\psi} &=& -\langle  \dot{M}_{0}, \psi\rangle  
            U_{0}\\
            Y_{\psi} (t=0) &=& 0.
        \end{array}
    \right .
    $$
\end{proposition}

\noindent {\bf (ii). Adjoint state technique}

    We have $\Lambda_h = tr_{t=T} \circ S$. Moreover, recall that 
    $$d\Lambda_h(0) : \left \{
        \begin{array}{ccc}
            C([0,T], W_h) &\longrightarrow
            &\m{F}(\interieur{\Omega_h})\\
            \varphi &\longrightarrow  & Y_{\varphi}(T)
        \end{array}
        \right .
        $$
     and $Y_{\varphi}$ is solution of $$
        \left \{
        \begin{array}{ccc}
            \partial_t Y_{\varphi} + M(0) y_{\varphi} &=& -\langle  \dot{M}(0), \varphi\rangle   U_0\\
            Y_{\varphi} (t=0) &=&0.
        \end{array}
        \right .
        $$
    Eventually, we denote $L : \varphi \longrightarrow Y_{\varphi}$. We clearly have $d\Lambda_h (0) = tr|_{t=T} \circ L$. Let us start by proving that the differential of $\Lambda_h$ at $0$ is surjective, using the adjoint state method.
     \begin{eqnarray}
        \nonumber d\Lambda_h (0) \textrm{ is surjective}
        &\Leftrightarrow & \{c \in \m{F}(\interieur{\Omega_h})^2;
        \forall \varphi \in \m{V} \langle  d\Lambda_h (0) \varphi, c \rangle   =0\}=
        \{0\}
    \end{eqnarray}
    and moreover,    
	\begin{align}
        \nonumber &\forall \varphi \in \m{V} & \langle  d\Lambda_h (0) \varphi, c\rangle=0\\
        \nonumber  \Leftrightarrow &\forall \varphi \in \m{V}&
        \langle  tr|_{t=T}(S(\varphi)),c\rangle  =0\\
\label{condition_surjective_ondes}
        \Leftrightarrow &\forall \varphi \in
        \m{V} & \langle  S(\varphi),tr|_{t=T}^*(c)\rangle   =0
     \end{align}

    Furthermore, we have $S(\varphi)= Y_{\varphi}$ is solution of the differential equation: 
    $$\partial_t Y_{\varphi} + M(0) y_{\varphi} = -M'_{\varphi} y_0$$

    \begin{definition}\label{def_etat_adjoint_ondes}
        The adjoint state $Y_{\varphi}$ and $tr|_{t=T}^*(c)$ is the unique solution $X$ of the equations:
        \begin{equation}\label{defadj_ondes}
        \left \{
            \begin{array}{ccc}
                -\partial_t X + M^* X &=& tr|_{t=T}^*(c)\\
                X(t=0)&=&0
            \end{array}
        \right .
        \end{equation}
    \end{definition}

With this definition, we replace \eqref{condition_surjective_ondes} $tr|_{t=T}^*(c)$ by its expression in function of the adjoint state $X$ \eqref{defadj_ondes}.\\
We deduce the following theorem:

    \begin{theorem}\label{surjection_theoreme_ondes}
        The differential $d\Lambda_h(0)$ of $\Lambda_h$ at $\varphi = 0$ is surjective if and only if we have the following uniqueness property:\\
        If $c\in \m{F}(\interieur{\Omega_h})$ is such that 
        \begin{equation}\label{condi_surj_theo_ondes}
        \langle  X, M'_{\varphi}y_0\rangle_{L^2(0,T;\m{F}(\interieur{\Omega_h}))} =
        0, \textrm{   } \forall \varphi \in \m{V}
        \end{equation}
        with 
        \begin{equation}
            \left \{
                \begin{array}{ccc}
                    -\partial_t X + M^*(0) X &=& tr|_{t=T}^*(c)\\
                    X(t=0)&=&0
                \end{array}
            \right .
        \end{equation}
        Then necessarily $c=0.$
    \end{theorem}

We now turn to compute the differential of $M$ at $0$
    \begin{proposition}\label{mprim}
    $$M'(0) =  \left (
        \begin{array}{cc}
            0 & 0\\
            A'(0) & 0
        \end{array}
    \right )$$
    Where $A'(0)$ is defined in \ref{aprim}.
    \end{proposition}
	A simple corollary of proposition~\ref{layer1} ensures that:
    \begin{proposition}\label{layer1_ondes}
        Assume that the reference state satisfies the discrete non-degeneracy condition on the whole boundary of $\Omega_h$. Then the relation \eqref{condi_surj_theo_ondes} implies that $$X|_{\Gamma^1} \equiv 0.$$
    \end{proposition}

    %
We are in a position to show the uniqueness property \eqref{condi_surj_theo_ondes}. From proposition \ref{layer1_ondes}, under the discrete non-degeneracy condition on the reference state $y_0$, we have  $X|_{\Gamma^1} \equiv 0$. We now show that this condition, together with the definition of the adjoint state \eqref{def_etat_adjoint_ondes}:
 $$ \left \{
            \begin{array}{ccc}
                -\partial_t X + M^*(0) X &=& tr|_{t=T}^*(c)\\
                X(t=0)&=&0
            \end{array}
        \right .
        $$

    \noindent implies that $X$ is identically vanishing, so $tr_{t=T}^*=0$ and $c=0$.\\ To this purpose, we will study the zeros propagation of $X$ on $\Omega_h$ from its boundaries.

    \begin{theorem}
        The relations:
        \begin{equation}
        \left \{
            \begin{array}{ccc}
                -\partial_t X + M^*(0) X &=& tr|_{t=T}^*(c)\\
                X(t=0)&=&0\\
                X|_{\Gamma^0} &=& 0\\
                X|_{\Gamma^1} &=& 0
            \end{array}
        \right .
        \end{equation}
        implies that  $X \equiv 0$ and $c = 0$.
    \end{theorem}
The proof uses exactly the same approach as the heat equation case. 

This results allows to prove the following controllability theorem. 

    \begin{theorem}
        Assume that the reference state $y_0$ defined in \eqref{eq:ReferenceState} satisfies the discrete non degeneracy condition \ref{NDD}. Let $Y_{\varphi}$ be the solution of the perturbed equation\eqref{eq:def-etat-pertube}. Let finally $Z_d = Y_0(T)$. Then there exist neighborhoods $\m{V}(0) \subset C([0,T];W_h)$ and $\m{V}(Z_d) \subset \fom^2$ such that $\forall Z \in \m{V}(Z_d)  \exists \varphi \in \m{V}(0)$ such that $Y_{\varphi}(T)=Z$.
    \end{theorem}

\section*{Conclusion}
In this paper we proved that the linearized heat equation was approximately controllable with respect to the shape of the domain, while the wave equation is locally exactly controllable. We addressed the same questions in the case of the semi-discrete equations in two dimensions in a square and we proved that the two types of equations are exactly controllable. Nevertheless, the methods we developed in this paper do not allow us to see the discrete control as an approximation of the continuous control in the wave equation. Another discretization method should be used to address this question, the mixed finite elements method. Indeed, we claim that one of the main obstacle to this interesting issue is the discretization method used, which does not behaves smoothly in the limit $h \to 0$. For instance we know (see \cite{NZ,IZ}) that in boundary control problem of the unidimensional wave equation, spurious modes with high frequency numerical oscillations appear and the observability constant tends to infinity when $h$ tends to $0$. It has been proved also that this semi-discrete model is not uniformly controllable in the limit $h \to 0$. 

Nevertheless, the results of Castro and Micu in \cite{CM} are promising. They studied a system based on a mixed finite element space semi-discretization the linear 1-D wave equation with a boundary control at one extreme. They show that the controls obtained with these semi-discrete systems can be chosen uniformly bounded in $L^2(0, T)$ and in such a way that they converge to the HUM control of the continuous wave equation, i.e. the minimal $L^2$-norm control. This result motivates to study in contrast to the classical finite element semi-discretization a mixed finite element scheme.

\noindent{\bf Acknowledgments:} The author warmly acknowledge Enrique Zuazua for hosting the author in his laboratory at the Universidad Aut{\'o}noma de Madrid, proposing the study and for insightful and helpful discussions. The author also warmly acknowledges Gr\'egoire Allaire for great and fruitful scientific discussions.


\begin{thebibliography}{99}

\bibitem{A}
\newblock G. Allaire,
\newblock ``Conception optimale de structures,''
\newblock Springer-Verlag, New York, 2007.

\bibitem{B}
\newblock Haim Brezis,
\newblock ``Analyse fonctionnelle,''
\newblock Masson, Paris, 1983.

\bibitem{CM}
\newblock C.~Castro and S.~Micu,
\newblock \emph{Boundary controllability of a linear semi-discrete 1-d wave equation
  derived from a mixed finite element method},
\newblock Numerische Mathematik, \textbf{102} (2006), 413--462.

\bibitem{C}
\newblock J.~C{\'e}a,
\newblock \emph{Numerical methods of optimum shape design,}
\newblock In E.~Haug and Jean C{\'e}a, editors,
\newblock ``Optimization of
  Distributed Parameter Structure'', Alphen aan den Rijn, Sijthoff and
  Noordhoff, the Netherlands, 1981.

\bibitem{cea:81}
\newblock J.~C{\'e}a.
\newblock \emph{Optimization of distributed parameter structures}.
\newblock NATO Advanced study Institutes series, 1981.

\bibitem{CZ}
\newblock D.~Chenais and E.~Zuazua,
\newblock \emph{Controlability of an elliptic equation and its finite difference
  approximation by the shape of the domain},
\newblock {Numerische Mathematik}, \textbf{95} (2003), 63--99.

\bibitem{chenais:75}
\newblock D. Chenais,
\newblock \emph{On the existence of a solution in a domain identification problem},
\newblock {Journal of Math. Analysis and Applications}, \textbf{52} (1975), 189--219.

\bibitem{Holmgren:1901}
\newblock E. Holmgren.
\newblock \emph{{\"U}ber systeme von linearen partiellen differentialgleichungen.}
\newblock {{\"O}fversigt af Kongl. Vetenskaps-Academien F{\"o}rhandlinger},
  \textbf{58} (1901), 91--103.

\bibitem{H}
\newblock L.~H{\"o}rmander,
\newblock ``Linear Partial Differential Operators'',
\newblock Springer-Verlag, New York, 1969.

\bibitem{IZ}
\newblock J.A. Infante and E. Zuazua.
\newblock \emph{Boundary observability for the space semi-discretizations of the 1 -
  d wave equation},
\newblock { Mathematical Modeling and Numerical Analysis}, \textbf{33} (1999), 407--438.

\bibitem{L1}
\newblock J.L. Lions.
\newblock ``Contr{\^o}le optimal de systemes gouvern{\'e}s par des {\'e}quations aux d{\'e}riv{\'e}es partielles'',
\newblock Dunod, Paris, 1968.

\bibitem{L}
\newblock J.L. Lions.
\newblock {``Contr{\^o}labilit{\'e} exacte, stabilisation et perturbations de
  syst{\`e}mes distribu{\'e}s (I and II)''}.
\newblock Masson, Paris, 1988.

\bibitem{L2}
\newblock J.L. Lions,
\newblock \emph{Exact controllability, stabilization and perturbations for
  distributed sytems},
\newblock {SIAM Review}, \textbf{30} (1988), 1--68.

\bibitem{Luenberger1969}
\newblock D.~G. Luenberger, 
\newblock ``{Optimization by Vector Space Methods}''.
\newblock John Wiley and Sons Inc., New York, 1969.

\bibitem{MS}
\newblock F.~Murat and J.~Simon.
\newblock \emph{Sur le Contr\^ole par un Domaine G\'eom\'etrique},
\newblock Publication du Laboratoire d'analyse num{\'e}rique, 189, Paris VI, 1976.

\bibitem{NZ}
\newblock M.~Negreanu and E.~Zuazua,
\newblock \emph{Uniform boundary controllability of a discrete 1-d wave equation,}
\newblock {Systems and Control Letters}, \textbf{48} (2003), 261--280.

\bibitem{rousselet:81}
\newblock B.~Rousselet,
\newblock \emph{Optimization of distributed parameter structures,}
\newblock {NATO advanced study institutes series}, \textbf{50} (1981), 1474--1501.

\bibitem{simon:80}
\newblock J.~Simon,
\newblock \emph{Differentiation with respect to the domain in boundary value
  problems},
\newblock {Numer. Func. Anal. Optim.} \textbf{2} (1980), 649--687.

\bibitem{S}
\newblock J.~Simon.
\newblock {\em Diferenciaci{\'o}n con respecto al dominio}.
\newblock Lecture notes, Universidad de Sevilla, 1989.

\bibitem{zolesio:81}
\newblock J.P. Zolesio.
\newblock \emph{Optimization of distributed parameter structures.}
\newblock {NATO advanced studies series} \textbf{50} (1981), 1089--1151.

\bibitem{Z1}
\newblock E.~Zuazua,
\newblock \emph{Some problems and results on the controllability of partial
  differential equations,}
\newblock In {\em European Congress of Mathematics}, Progress in Mathematics.
  Birkh\"auser, 1996.

\bibitem{Z6}
\newblock E.~Zuazua,
\newblock \emph{Some new results related to the null controllability of the 1-d heat
  equation},
\newblock {\em Seminaire X-EDP, Ecole Polytechnique,}, \textbf{VIII} (1997--1998), 1--22.

\bibitem{Z4}
\newblock E.~Zuazua,
\newblock \emph{Boundary observability for the finite-difference space
  semi-discretizations of the 2-d wave equation in the square,}
\newblock {J. Math. Pures Appl.} \textbf{78} (1999), 523--563.

\bibitem{Z3}
\newblock E.~Zuazua,
\newblock \emph{Observability of the 1-d waves in heterogenous and semi discrete
  media},
\newblock In {\em Advances in Structural Control}. CIMNE, 1999.

\bibitem{Z2}
\newblock E.~Zuazua,
\newblock \emph{Controllability of partial differential equation and its
  semi-discrete approximations},
\newblock {Discrete and Continuous Dynamical Systems} \textbf{8} (2002), 469--513.

\end{thebibliography}

\medskip
Received xxxx 20xx; revised xxxx 20xx.
\medskip
\end{document}